\documentclass[11pt]{article}

\def\mathbb{\Bbb}
\usepackage{amsfonts,latexsym,amstext,theorem,euscript,array,enumerate,mathrsfs}
\usepackage{amsmath}
\usepackage{amssymb}

\usepackage{comment}

\newtheorem{Theorem}{Theorem}[part]
\newtheorem{Definition}{Definition}[part]
\newtheorem{Proposition}{Proposition}[part]

\newtheorem{Lemma}{Lemma}[part]

\newtheorem{Remark}{Remark}[part]


\def\qed{{\hfill\hbox{\enspace${ \square}$}} \smallskip}
\def\sqr#1#2{{\vcenter{\vbox{\hrule height .#2pt \hbox{\vrule
 width .#2pt height#1pt \kern#1pt \vrule
width .#2pt} \hrule height .#2pt}}}}
\def\square{\mathchoice\sqr54\sqr54\sqr{4.1}3\sqr{3.5}3}

\def\ds{\begin{displaystyle}}
\def\eds{\end{displaystyle}}
\def\dis{\displaystyle }
\def\<{\langle }
\def\>{\rangle }

\def\R{\mathbb R}

\def\E{\mathbb E}
\def\P{\mathbb P}
\def\Q{\mathbb Q}
\def\F{\mathbb F}

\def\D{\mathbb D}

\def\H{\mathbb H}
\def\G{\mathbb G}
\def\A{\mathbb A}

\def\cala{{\cal A}}
\def\calb{{\cal B}}
\def\calc{{\cal C}}

\def\calf{{\cal F}}
\def\calg{{\cal G}}
\def\calh{{\cal H}}

\def\calk{{\cal K}}
\def\calm{{\cal M}}

\def\calp{{\cal P}}

\def\calv{{\cal V}}

\def\cals{{\cal S}}
\def\bfB{{\bf B}}
\def\bfC{{\bf C}}

\def\beqs{\begin{eqnarray*}}
\def\enqs{\end{eqnarray*}}
\def\beq{\begin{eqnarray}}
\def\enq{\end{eqnarray}}

\title{Dual and backward SDE representation for \\ optimal control
of  non-Markovian   SDEs}

\date{}

\author{Marco Fuhrman
\\Politecnico di Milano,
Dipartimento di Matematica\\
via Bonardi 9, 20133 Milano, Italy\\
marco.fuhrman@polimi.it
\\
\\
Huy\^{e}n Pham
\\
LPMA - Universit\'{e} Paris Diderot \\
Batiment Sophie Germain, Case 7012 \\
13 rue Albert Einstein, 75205 Paris Cedex 13\\
and CREST-ENSAE \\
pham@math.univ-paris-diderot.fr
}

\begin{document}

\maketitle

\begin{abstract}
We study optimal stochastic control problem for  non-Markovian  stochastic differential equations (SDEs) where the drift, diffusion coefficients, and gain functionals are 
path-dependent, and importantly we do not make any ellipticity assumption on the SDE.   
We develop a controls randomization approach, and prove that the value function can be reformulated under a family of  dominated measures
on an enlarged filtered probability space.  This value function is then characterized by a backward SDE with nonpositive jumps under a single probability measure, which can be viewed as 
a path-dependent version of the Hamilton-Jacobi-Bellman equation, and an extension to $G$-expectation. 
\end{abstract}

\vspace{5mm}


\vspace{5mm}

\noindent {\bf MSC Classification (2010):} 60H10, 93E20.

\newpage

\section{Introduction}

We consider non-Markovian controlled stochastic differential equations (SDEs) of the form
\beq \label{SDEpath}
dX_s &=&  b_s(X,\alpha_s)\,ds  +     \sigma_s(X,\alpha_s)\,dW_s, \;\;\; 0 \leq s \leq T,
\enq
where $W$ is a  $n$-dimensional Wiener process, $\alpha$ is a progressive control process, and the drift and diffusion coefficients
$b$ and $\sigma$ may depend on the trajectory of the solution $X$ valued in $\R^d$ in a non-anticipative way.  Given initial conditions determined in our  context by $t$ $\in$ $[0,T]$, and $x$ $\in$ $\bfC^d$, the set of
continuous function from $[0,T]$ into $\R^d$,  we denote by $X^{t,x,\alpha}$  the solution to  \eqref{SDEpath} associated to the control $\alpha$, and starting from $X_s$ $=$ $x(s)$ for $s$ $\in$
$[0,t]$.  We are then interested in the value function for the optimal stochastic control problem:
\beq \label{defvaluenonmarkov}
v(t,x) &=& \sup_\alpha \E \Big[ \int_t^T f_s(X^{t,x,\alpha},\alpha_s) ds + g(X^{t,x,\alpha}) \Big],
\enq
where the running and terminal reward functionals $f$ and  $g$  may also depend on the past trajectory of  the solution $X$.

In the Markovian framework (see e.g. \cite{kr}), i.e. when $b_s$, $\sigma_s$, $f_s$ depend on $X$ only through its current value $X_s$,  and $g$ only on $X_T$,
the value function also depends at time $t$  only on  the current state value $X_t$ $=$ $x(t)$, hence is a  deterministic function on $[0,T]\times\R^d$.  By the dynamic programming approach, the value function is then described by  the Hamilton-Jacobi-Bellman (HJB)  partial differential equation (PDE), which is satisfied in general in the viscosity sense, and which characterizes the control problem once we have a uniqueness result for the HJB PDE.  We refer to the monographs \cite{fleson94} or \cite{pha09}  for a detailed exposition of this theory of dynamic programming  and viscosity solutions for stochastic optimal control.

The representation of stochastic control problem and HJB equation  has been also  developing by means of backward stochastic differential equations (BSDEs). In the Markovian case where the controller can affect only the drift coefficient,  the HJB equation is a semi-linear PDE and is known to be related to a standard BSDE, see  \cite{parpen92}.  The controlled diffusion case, arising typically in finance in uncertain volatility models,  leads to a fully nonlinear HJB PDE, and can be represented by  second-order BSDE (2BSDE), as introduced in \cite{cheetal07} and  \cite{sontouzha11},  whose  basic idea is to require that the solution verifies the equation  almost surely for every probability measure in a non dominated class of mutually singular measures. This theory is closely related to the notion of nonlinear and $G$-expectation, see \cite{pen06}, but requires a nondegeneracy condition on the diffusion coefficient together with some constraint between drift and diffusion.  The general case without any ellipticity assumption on the controlled diffusion is addressed in \cite{khph}, where it is proved that 
fully nonlinear HJB equations can be  represented by a class of BSDE with nonpositive jumps. The basic idea, following \cite{khaetal10} (see also \cite{bou09} for optimal switching problem),  is to randomize the control process $\alpha$ by replacing it by an uncontrolled pure jump process associated to a Poisson random measure, and then to  constrain  the jumps-component solution to the BSDE driven by  Brownian motion and  Poisson random  measure,  to remain nonpositive, by adding a nondecreasing process in a minimal way.  A key  feature of this class of BSDEs is its formulation under a single probability measure like for standard BSDE  in contrast with 2BSDEs, thus avoiding technical issues in quasi-sure analysis. It is then proved in \cite{khph} that the minimal solution to the BSDE with nonpositive jumps satisfies the nonlinear HJB equation, so that it coincides with the value function of the Markovian stochastic control problem, once one has at disposal a  uniqueness result for this HJB PDE (see \cite{craishlio92} for a review on comparison results for viscosity solutions to nonlinear  PDEs).

 The main goal of this paper is to extend the result of \cite{khph} to the non-Markovian framework.  More precisely, we aim to prove that the value function in \eqref{defvaluenonmarkov} may be represented in terms of a  BSDE with nonpositive jumps, which can then be seen as a non-Markovian version of  HJB equation.    We use a controls randomization approach, by replacing the control process $\alpha$ by a pure-jump process associated to a Poisson random measure independent of the Wiener process, with fixed finite intensity measure.  We then  show that  the value function in its weak formulation can be written equivalently as a dual control problem under a family of
dominated (and even equivalent) probability measures on an enlarged probability space, whose effect is to change the intensity measure of the Poisson random measure.  By means of this dual representation, we are finally able to relate the minimal solution to the BSDE with nonpositive jumps to the original value function \eqref{defvaluenonmarkov}.  The arguments in this paper for proving this connection are quite different from the Markovian case studied in \cite{khph}. Indeed, this connection is shown in that paper through the HJB equation, which is satisfied both by the value function and the minimal solution, and thus requires a uniqueness result.
Here, we prove this connection through the dual control problem by purely probabilistic arguments. The main issue is to  approximate  continuous control processes by pure jump processes associated to random measures with compensator absolutely continuous with respect to a given finite intensity measure.
In particular, we  do not rely on the path-dependent HJB equation associated by dynamic programming principle to  the value function in the non-Markovian context,  thus circumventing delicate issues of dynamic programming (as originally studied in \cite{elk79} for general non-Markovian stochastic control problems),  viscosity solutions and  comparison principles for fully nonlinear path-dependent PDEs,
as recently stu\-died in   \cite{pen11}, \cite{ekretal11} and  \cite{tanzha12}, see also \cite{FGSbook} for HJB equations in infinite dimension arising typically for  stochastic systems with delays. 
This suggests in particular an original approach to derive the HJB equation for value function of stochastic control problem from the BSDE representation,  hence without dynamic programming principle.

 We mention  that optimal control for path-dependent SDEs  was also recently studied in \cite{nut12} by adopting a quasi-sure formulation approach, which allows the author to prove a pathwise dynamic programming principle, and to derive a 2BSDE satisfied by the value function.  However,  the results are obtained essentially under a non degeneracy condition on the matrix  diffusion coefficient and when control cannot affect independently  drift and diffusion, see for details Assumption 2.1 and Remark 2.2 in \cite{nut12}.  Our results do not require any non degeneracy condition on $\sigma$,  and include the case of control both on drift and diffusion coefficient arising for instance in portfolio optimization problems.

 The rest of the  paper is organized as follows.  In Section 2,  we detail the controlled path-dependent SDE, and introduce the corresponding value function in its weak formulation.
 Section 3 presents  the main results of the paper. We  formulate the value function by means of a dual control problem over changes of  dominated measures on an enlarged probability space.
 This dual representation  allows us to characterize the value function as the solution to a BSDE with nonpositive jumps.  The proofs  are  reported in Section 4.
 Finally, we collect in Appendix some useful results about random measures and their compensators.

\section{Control of path-dependent SDEs}

\setcounter{equation}{0} \setcounter{Assumption}{0}
\setcounter{Theorem}{0} \setcounter{Proposition}{0}
\setcounter{Corollary}{0} \setcounter{Lemma}{0}
\setcounter{Definition}{0} \setcounter{Remark}{0}

We introduce in this section the path-dependent control setting and assumptions, and we define the value function for the associated optimal control problem.

\subsection{Non-markovian controlled SDE}

Let $A$, the control space, be a  Lusin space (some authors call it a Borel space), i.e. a topological space  homeomorphic to a  Borel subset of a Polish space,  endowed with a  metric, denoted by $\rho$. We may assume without loss of generality that $\rho(a,a')$ $<$ $1$, for any $a,a'$ $\in$ $A$, by replacing otherwise the initial metric by the equivalent one: $\rho/(1+\rho)$.  We denote by $\calb(A)$ the Borel $\sigma$-algebra of $A$. In the  sequel,  we also need to consider the set of all positive finite measures on $(A,\calb(A))$ with full topological support, which will be denoted by $\calm_f(A)$.  We note that $A$ can be a finite or countable set, a Borel subset of $\R^q$, or more generally
any Polish space.
%(i.e. separable topological space with the topology induced by a complete metric).

In order to specify the measurability assumptions on the coefficients of the controlled equation, and of the control problem over a fixed horizon $T$ $<$
$\infty$,  we introduce  the path space $\bfC^d$  of continuous maps from $[0,T]$ to $\R^d$   and we equip $\bfC^d$ with the usual supremum norm
$\|x\|_{_\infty}$ $=$ $x^*_T$, where we set  $x^*_t$ $:=$ $\sup_{s\in [0,t]}|x(s)|$, for $t$ $\in$ $[0,T]$, and $x$ $\in$ $\bfC^d$. We  define the filtration
$(\calc_t)_{t\in[ 0,T]}$, where $\calc_t$ is the $\sigma$-algebra generated by the canonical coordinate maps $\bfC^d\to \R^d$,
$x(\cdot)\mapsto x(s)$ up to time $t$:
\beqs
\calc_t  &:=  &  \sigma \{ x(\cdot)\mapsto x(s)\; :\, s\in [0,t] \}.
\enqs
Let $Prog(\bfC^d)$ denote the progressive $\sigma$-algebra in $[0,T]\times\bfC^d$ with respect to $(\calc_t)$.

The drift and diffusion coefficients
\begin{center}
\begin{tabular}{ccccccc}
$[0,T]\times\bfC^d\times A $ & $\longrightarrow$ &  $\R^d$, &  &  $[0,T]\times\bfC^d\times A$ & $\longrightarrow$ & $\R^{d\times n}$ \\
$(t,x,a)$ &$\longmapsto$& $b_t(x,a)$, & & $(t,x,a)$ &$\longmapsto$ & $\sigma_t(x,a)$
\end{tabular}
\end{center}
are $Prog(\bfC^d)\otimes \calb(A)$-measurable.  This measurability requirement is a standard one for ensuring that the stochastic
differential equations (SDEs) with path-dependent coefficients is well-posed.  We shall make the usual assumption:

\vspace{3mm}

\noindent {\bf (H1)}
\begin{itemize}
\item [(i)]  For all $t$ $\in$ $[0,T]$, and $x$ $\in$ $\bfC^d$, the functions $b_t(x,a)$, and $\sigma_t(x,a)$ are continuous in $A$.
\item [(ii)] There exists a nonnegative constant  $K_1$ such that
\beq
|b_t(x,a) - b_t(x',a)| + |\sigma_t(x,a)-\sigma_t(x',a)| & \leq & K _1 (x-x')^*_t, \label{lipbsig} \\
|b_t(0,a)| + |\sigma_t(0,a)| & \leq & K_1,  \label{borbsig}
\enq
for all $(t,x,x',a)$ $\in$ $[0,T]\times\bfC^d\times\bfC^d\times A$.
\end{itemize}

\vspace{2mm}

We now formulate the controlled path-dependent SDE.  Borrowing some terminology from \cite{rowi},
by an {\em admissible set-up} (or simply a set-up) we mean
\beqs
\A &=& (\Omega, \calf, \G, \Q, W),
\enqs
where  $(\Omega, \calf,  \Q)$ is a probability space equipped with a filtration  $\G$ $=$ $(\calg_t)_{t\ge 0}$ satisfying the {\it usual conditions},
and $W$ $=$ $(W_t)_{0\leq t\leq T}$ is an $n$-dimensional standard $(\Q,\G)$-Wiener process.
Notice that $\G$ is not necessarily the natural filtration of $W$.
%this simplifies some of the arguments below.
We define the space of $\A$-admissible controls, denoted $\cala(\A)$,   as the set  of processes  defined on $[0,T]\times\Omega$,
valued in  $A$, which are  progressively measurable  (for short, progressive) with respect to  $\G$. For fixed $t$ $\in$ $[0,T]$, $x$ $\in$ $\bfC^d$,
and given $\alpha$ $\in$ $\cala(\A)$, we consider the stochastic differential equation:
\begin{equation}\label{stateeq}
    \left\{\begin{array}{lcl}
X_s \; = \; x(s),  &   & s\in[0,t], \\\dis
X_s \; = \; x(t) +  \int_t^s b_u(X,\alpha_u)\,du  + \int_t^s       \sigma_u(X,\alpha_u)\,dW_u, &  &  s \in [t,T].
\end{array}\right.
\end{equation}
By standard results (see e.g. \cite{rowi}, Thm V. 11.2),
under {\bf (H1)},  there exists a unique $\G$-adapted strong solution $X$ $=$ $(X_s)_{0\leq s \leq T}$   to \eqref{stateeq}
with continuous trajectories and satisfying, for every $p\in [1,\infty)$,
\beq\label{solusp}
    \E^\Q\,\Big[\sup_{s\in [t,T]}|X_s|^p\Big]  & \leq &  C \Big(1+ (x^*_t )^p\Big),
\enq
(where $\E^\Q$ denotes of course the expectation under $\Q$)
for some constant $C$ depending only on $p,T$, and $K_1$  as defined in {\bf (H1)}.
We will denote the solution $X^{t,x,\A,\alpha}$ to stress dependence on these parameters. Notice that $\{X_s^{t,x,\A,\alpha}$, $t\leq s\leq T\}$, depends on $x$ only on the past trajectory $\{ x(u), 0 \leq u \leq t\}$.

\begin{Remark}
\begin{em} \label{deterministiccase}
It is worth  to mention that in contrast with \cite{sontouzha11}, \cite{nut12},  no non-degeneracy assumption on the diffusion coefficient $\sigma$, nor  specific condition between the drift and the diffusion coefficient  is imposed.
In particular, we may control independently both drift and diffusion, and it may well happen that some lines or columns of $\sigma$ are equal to zero, and even  $\sigma=0$, in which case we have results   for   deterministic control problems.
\end{em}
\end{Remark}

\begin{Remark}\begin{em}
\label{winthecoefficients}
One may consider a priori more general  non-Markovian controlled equations of the form
\beqs
dX_t &=&  b_t(X,W,\alpha_t)\,dt  +     \sigma_t(X,W,\alpha_t)\,dW_t,
\enqs
i.e. where the drift and diffusion coefficients are non-anticipative functionals
of the trajectory of the Wiener process $W$ as well. However, this generality
is only apparent, since we can adopt the following standard procedure to
reduce to the case presented above: we introduce a second state component
$Y=W$ and consider the equivalent controlled system
$$
\left\{\begin{array}{l}
dX_t \; = \;  b_t(X,Y,\alpha_t)\,dt  +     \sigma_t(X,Y,\alpha_t)\,dW_t, \\
dY_t \; = \; dW_t
\end{array}
\right.
$$
which is of the form considered above, but with an enlarged state $(X,Y)$. It is easy to formulate assumptions on
$b_t(X,W,\alpha_t)$, $\sigma_t(X,W,\alpha_t)$ that allow to verify the requirements in \eqref{lipbsig}-\eqref{borbsig}
on the resulting controlled system, so we omit the details.  We only remark that the fact that the latter has degenerate noise does not prevent the possibility of applying our results, as noted above.
\qed
\end{em}
\end{Remark}

\subsection{The value function}

We are given a running and terminal reward function:
\begin{center}
\begin{tabular}{ccccccc}
$[0,T]\times\bfC^d\times A $ & $\longrightarrow$ &  $\R$, &  &  $\bfC^d$ & $\longrightarrow$ & $\R$ \\
$(t,x,a)$ &$\longmapsto$& $f_t(x,a)$, & & $x$ &$\longmapsto$ & $g(x)$,
\end{tabular}
\end{center}
which are respectively $Prog(\bfC^d)\otimes \calb(A)$-measurable, and $\calc_T$-measurable, and we assume:

\vspace{3mm}

\noindent {\bf (H2)}
\begin{itemize}
\item [(i)]  For all $t$ $\in$ $[0,T]$,  the function $f_t(x,a)$ is continuous in $(x,a)$ $\in$ $\bfC^d\times A$, and the function $g$ is continuous on $\bfC^d$,  \item [(ii)] There exist  nonnegative  constants  $K_2$ and $m$ such that
\beqs
|f_t(x,a)| + |g(x)| & \leq & K_2\big(1 + \|x\|_{_\infty}^m \big),
\enqs
for all $(t,x,a)$ $\in$ $[0,T]\times\bfC^d\times A$.
\end{itemize}

\vspace{2mm}

 We then define the gain functional:
 \beqs
 J(t,x,\A,\alpha) &=& \E^\Q\, \Big[\int_t^T f_s(X^{t,x,\A,\alpha},\alpha_s)\,ds +  g(X^{t,x,\A,\alpha})\Big],
 \enqs
 for $(t,x)$ $\in$ $[0,T]\times\bfC^d$, $\A$ $=$ $(\Omega, \calf,\G, \Q, W)$ a set-up, $\alpha$ $\in$ $\cala(\A)$, and the value function in its weak formulation
 as the supremum over all admissible set-ups and controls:
\beq \label{defvalue}
v(t,x) &=&  \sup_\A \sup_{\alpha\in\cala(\A)} J(s,x,\A,\alpha), \;\;\; (t,x) \in [0,T]\times\bfC^d.
\enq
Due to \eqref{solusp} and the polynomial growth condition on $f,g$ in {\bf (H2)},  it is easy to check that $v$ is always finite, and satisfies actually:
\beq \label{growthv}
|v(t,x)| & \leq & K(1+  |x^*_t|^m), \;\;\;  (t,x) \in [0,T]\times\bfC^d,
\enq
for some positive constant $K$.  Thus, $v$ is a real-valued  function defined on $[0,T]\times\bfC^d$. Moreover,
since $\{X_s^{t,x,\A,\alpha}$, $t\leq s\leq T\}$, depends on $x$ only on the past trajectory $\{ x(u), 0 \leq u \leq t\}$, this is also the case for
$v(t,x)$.  In other words, $v(t,.)$ is $\calc_t$-measurable for all $t$ $\in$ $[0,T]$.  We do not address here the question of joint measurability of $v$ in its arguments, i.e. whether $v$ is
$Prog(\bfC^d)$-measurable, since it is not needed for our purpose.
%\marginpar{Other comments about joint measurability?}
We simply mention that  this issue is already not trivial in the Markovian case, and usually relies  on  a measurable selection theorem.

\vspace{2mm}

\begin{Remark}
{\rm One could also consider the optimal control problem in the strong formulation,
i.e. the search for an optimal control $\alpha\in\cala(\A)$ in a given set-up $\A$,
and the corresponding value function (in general depending on $\A$) which
is defined as in
\eqref{defvalue} but dropping the supremum with respect to $\A$:  see e.g.
\cite{FGSbook} Section IV.2,  or \cite{yozh} Sections 2.4.1-2.4.2
for  detailed formulations. In the Markovian framework, when a verification
theorem for the HJB equation holds under appropriate conditions,
the value functions for the weak and strong formulations are known to be
the same: see \cite{FGSbook} Remark IV.3.2. In this paper we will not address
similar issues for the present path-dependent context: instead,
we are rather interested in relating
the function $v(t,x)$, as defined above, to a dual control problem and
to a suitable BSDE.
\qed
}
\end{Remark}

\begin{Remark}
{\rm Fix an initial condition $t$ $=$ $0$, $x$ $=$ $0$,  a set-up $\A$, and given a control $\alpha$ $\in$ $\cala(\A)$,
denote by $\Q^\alpha(\A)$ the distribution  of $X^{0,0,\A,\alpha}$, which can be seen as a probability measure on the canonical space $\bfC^d$.
Then, the family  $\{\Q^\alpha(\A), \alpha \in \cala(\A)\}$  is  not dominated in general when the diffusion coefficient depends on the control $\alpha$. In particular, when $b$ $=$ $0$, $\sigma_t(x,a)$ $=$ $a$, so that the SDE \eqref{stateeq} degenerates to  a stochastic integral, and for $f$ $=$ $0$, we see that the value function \eqref{defvalue} falls into the class of sublinear expectations studied in \cite{denetal}. More precisely, by considering $g$ as a random variable on the canonical space, we have:
\beqs
v(0,0) &=& \sup_\A \sup_{\alpha\in\cala(\A)} \E^{\Q^\alpha(\A)}[g],
\enqs
so that the mapping $g$ $\mapsto$ $v(0,0)$ may be viewed as a generalization of $G$-expectation \cite{pen06},
where the  volatility  $\alpha_t$ of the canonical process is uncertain, valued in $A$.
%Our framework slightly differs from \cite{pen06}
%in the sense that we allow the uncertain volatility $\alpha$  to be progressive with respect to
%a filtration not necessarily equal to the filtration generated by the canonical process.
\qed
}
\end{Remark}

\section{Dual control problem and BSDE representation}

\setcounter{equation}{0} \setcounter{Assumption}{0}
\setcounter{Theorem}{0} \setcounter{Proposition}{0}
\setcounter{Corollary}{0} \setcounter{Lemma}{0}
\setcounter{Definition}{0} \setcounter{Remark}{0}

In this section, we provide a dual representation of the control problem \eqref{defvalue}  by randomization of the controls.  This will allow us to
characterize the value function as the solution to  a backward stochastic differential equation (BSDE) with nonpositive jumps, formulated under a single probability measure.
This should be understood as a non-Markovian analog of
fully nonlinear Hamilton-Jacobi-Bellman equation, and an alternative to 2BSDE (\cite{sontouzha11}, \cite{nut12}) under more general conditions, see Remark \ref{deterministiccase}.  We shall assume
that {\bf (H1)} and {\bf (H2)} hold throughout this section.
%or viscosity solutions of path-dependent PDEs as recently studied in \cite{pen11}, \cite{ekretal11}.

\subsection{Dual representation with dominated measures} \label{paradual}

To state the dual control problem we initially fix a finite measure $\lambda$  on $(A,\calb(A))$ with full topological support,
i.e. an element of   $\calm_f(A)$, according to our previous notation.
By a {\em dual admissible set-up} (or simply a dual set-up) we mean
\beqs
\D_\lambda &=& (\Omega, \calf,  \P, W,\mu_\lambda),
\enqs
 where  $(\Omega, \calf, \P)$ is a probability space,    $W$ $=$ $(W_t)_{0\leq t\leq T}$  is an $n$-dimensional standard Wiener process,
$\mu_\lambda=\sum_{n\ge1}\delta_{(S_n,\eta_n)}$
is a  Poisson random measure in $[0,\infty)\times A$ with compensator $\lambda(da)dt$,  independent of $W$.
We denote by  $\F$ $=$ $(\calf_t)_{t\ge 0}$ the completion of the  natural filtration of the pair $(W,\mu_\lambda)$.
Although it depends on $\D_\lambda$, we do not make it explicit in the notation.

Let $\calv(\D_\lambda)$ be the class of random fields $\nu_t(\omega,a)$ $:$ $[0,\infty)\times\Omega\times A\to (0,\infty)$
which are $\calp(\F )\otimes \calb(A)$-measurable, where $\calp(\F)$ denotes the predictable $\sigma$-algebra associated to $\F$, and essentially bounded with respect to the measure
$dt\otimes d\P\otimes \lambda(da)$.
For $\nu\in \calv(\D_\lambda)$, the Dol\'eans-Dade exponential process
\beq \label{LGirsanov}
   L^\nu_t & := & \exp\Big(\int_0^t\int_A(1-\nu_s(a))\, \lambda(da)\,ds\Big) \prod_{S_n\le t} \nu_{_{S_n}}(\eta_n),\qquad t\ge0,
\enq
is a positive $(\P,\F)$-martingale on $[0,T]$, i.e. $\E[L_T^\nu]$ $=$ $1$  (since $\nu$ is bounded and $\lambda(A)$ $<$ $\infty$),  and defines a probability measure $\P^\nu$  on $(\Omega,\calf)$,
equivalent to $\P$, by setting
$\P^\nu(d\omega)$ $=$ $L_T^\nu(\omega)\P(d\omega)$.
%\beqs
%\left. \frac{d\P^\nu}{d\P} \right|_{\calf_t} &=& L_t^\nu, \;\;\; t \geq 0.
%\enqs
We recall that,  by Girsanov's theorem (see e.g. \cite{ja} Theorem 4.5),
under $\P^\nu$ the random measure $\mu_\lambda$ admits compensator $\nu_t(a)\,\lambda(da)dt$ on $[0,T]\times A$, and $W$ remains a Wiener process
independent of $\mu_\lambda$.

In order to define the gain functional, in addition to $\lambda\in   \calm_f(A)$ we also fix  an element $a\in A$, a starting time $t\in [0,T]$  and an initial condition $x$ $\in$ $\bfC^d$. Then we define
 a pair of processes $(I,X)$ $=$ $(I_s,X_s)_{0\leq s \leq T}$ as the unique strong solution to
\begin{equation}\label{defxi}
\left\{
    \begin{array}{llll}
I_s &=&  a, & s\in [0,t], \\
I_s &=& \dis a + \int_{(t,s]} \int_A (a'-I_{u-})\,\mu_\lambda(du,da'), & s \in [t,T], \\
X_s &=& x(s),& s \in[0,t], \\
X_s &=& \dis x(t) + \int_t^s b_u(X,I_u)\,du  +\int_t^s  \sigma_u(X,I_u)\,dW_u,&  s\in [t,T].
\end{array}
\right.
\end{equation}
We note that $I$ is determined by the initial point $a$ and the restriction of $\mu_\lambda$ to $ (t,T]\times A$: more precisely, letting $N_t$ $=$ $\sum_{n\ge1}1_{S_n\le t}$
denote the number of jumps of $\mu_\lambda$ in the time interval $[0,t]$, we have the explicit formula
\beqs
I_s &=& a\,1_{[0,S_{N_s+1} )}(s) + \sum_{n\ge N_s+1}\eta_n\,1_{[S_n,S_{n+1})}(s), \qquad s\in [0,T].
\enqs
The uncontrolled pure-jump process $I$ valued in $A$,  should be understood as a randomization of the control $\alpha$ in the primal problem defined in the previous paragraphs.
Instead of $(I,X)$, we may  write $(I^{t,a,\lambda,\D_\lambda},X^{t,x,a,\lambda,\D_\lambda})$ to stress dependence on these parameters.

We then  introduce the dual gain functional
\beqs
J^*(t,x,a,\lambda,\D_\lambda, \nu) &=& \E^\nu \Big[ \int_t^T f_s(X^{t,x,a,\lambda,\D_\lambda},I_s^{t,a,\lambda,\D_\lambda}) \,ds + g(X^{t,x,a,\lambda,\D_\lambda})\Big],
\enqs
(here $\E^{\nu} $ denotes the expectation under $\P^{\nu}$) for $(t,x,a)$ $\in$ $[0,T]\times\bfC^d\times A$, $\lambda\in \calm_f(A)$, $\D_\lambda$ a dual set-up, and $\nu$ $\in$ $\calv(\D_\lambda)$.
The intuitive  interpretation is that in the dual control problem we first replace the control $\alpha$ by a Poisson point process $I$, and we can then
 control the system by modifying the intensity of $I$. More precisely, we are able to  change its compensator $\lambda(da)dt$  into $\nu_t(a)\,\lambda(da)dt$,  by choosing the density
$\nu$ in the class of bounded positive and predictable random fields.

We finally introduce two  dual value functions. The first one is defined by
 \beq \label{defdual1}
 v_1^*(t,x,a,\lambda, \D_\lambda) &=& \sup_{\nu\in\calv(\D_\lambda)} J^*(t,x,a,\lambda,\D_\lambda, \nu),
\enq
and corresponds to optimizing with respect to every choice of $\nu$ in the fixed set-up $\D_\lambda$ and for a fixed  $\lambda\in \calm_f(A)$, i.e. it corresponds to the strong  formulation of the dual optimal control problem.
The second one is the value function in the weak formulation, where the set-up is part of the control:
\beq \label{defdual2}
v_2^*(t,x,a,\lambda) &=& \sup_{\D_\lambda}  v_1^*(t,x,a,\lambda, \D_\lambda).
% \sup_{\nu\in\calv(\D_\lambda)} J^*(t,x,a,\lambda,\D_\lambda, \nu).
\enq
We shall see later that the strong and weak dual formulations are equivalent, i.e. $v_1^*$ does not depend on $\D_\lambda$, and so is equal to $v_2^*$.  For the moment, we can easily check, as in the original control problem, that $v_1^*$ and $v_2^*$ are always finite, satisfying actually the same growth condition
\eqref{growthv}, and $v_1^*(t,x,.)$ and $v_2^*(t,x,.)$ depend on $x$ $\in$ $\bfC^d$  only via the past trajectory $\{x(s), 0\leq s \leq t\}$.
%Thus $v_d$ is a real function defined for $u\in [0,T]$, $x(\cdot)\in\bfC$, $a\in A$, $\lambda\in \calm_f(A)$,
%while $\tilde v_d$ is a function also of the dual set-up $\D_\lambda$.

\vspace{2mm}

Our first main result is to connect the primal control problem to the weak  dual one.

\begin{Theorem}\label{main}
%Let Assumptions {\bf (H1)} and {\bf (H2)} hold.
We have
 \beqs
  v(t,x) &=&    v_2^*(t,x,a,\lambda), \;\;\;\;\; (t,x) \in [0,T]\times\bfC^d,
\enqs
for all $a\in A$, $\lambda\in\calm_f(A)$.  In particular, $v_2^*(t,x,a,\lambda)$ does not depend on $(a,\lambda)$.
\end{Theorem}

\begin{Remark}
{\rm The dual problem is a control problem over equivalent probability measures $\P^\nu$, whose effect is to change the intensity measure of the pure-jump component $I$.
Theorem \ref{main} then formally means that one can formulate  the primal control problem originally written  in a non dominated Wiener space framework into a dominated framework  by enlarging
the filtered probability space with an additional Poisson random measure. Moreover, the result is invariant with respect to the choice of the intensity measure  for the Poisson random measure. 
\qed
}
\end{Remark}

\subsection{BSDE characterization}

Throughout this section,  we fix $\lambda\in\calm_f(A)$ and the initial conditions  $t$ $\in$ $[0,T]$,   $x$ $\in$ $\bfC^d$ and $a\in A$.
We consider a dual  set-up $\D_\lambda=(\Omega, \calf,  \P, W,\mu_\lambda)$,
denote by $\F$ the completion of the natural filtration of $(W,\mu_\lambda)$, and
define  a pair of processes $(I,X)$  as the solution to the system  \eqref{defxi}, dropping their dependence on  $t,x,a,\lambda,\D_\lambda$ in the notation.

Following \cite{khph}, we say that a quadruple $(Y,Z,U,K)$ is a solution to the  BSDE with
nonpo\-si\-ti\-ve jumps:
\begin{equation} \label{BSDEnonpos}
\left\{
\begin{array}{ccl}
Y_s &=& \dis g(X) + \int_s^T f_r(X,I_r) dr  + K_T - K_s \\
       &   &\dis -\int_s^T Z_r dW_r -  \int_s^T \int_A  U_r(a) \mu_\lambda(dr,da),  \;\;\; s \in [0,T],  \\
U_s(a) & \leq & 0,
\end{array}
\right.
\end{equation}
if $Y$ $\in$ $\cals^2$,  the space of  c\`ad-l\`ag $\F$-adapted processes  $Y$ satisfying $\|Y\|^2$ $:=$ $\E[\sup_{s\in[0,T]}|Y_s|^2]$ $<$ $\infty$,  identified up to indistinguishability, $Z$ $\in$  $L^2(W)$,  the space of $\F$-predictable processes with values in $\R^n$ such that $\|Z\|^2_{L^2(W)}$ $:=$ $\E[\int_0^T|Z_s|^2ds]$ $<$ $\infty$,  identified up to $ds\otimes d\P$-a.e. equality,
$U$ $\in$ $L^2(\tilde\mu)$,  the space of  $\calp(\F)\otimes \calb(A)$-measurable real-valued  processes $U$ such that
$\|U\|^2_{L^2(\tilde\mu)}$ $:=$ $\E[\int_0^T\int_A|U_s(a)|^2\lambda(da)\,ds]$ $<$ $\infty$,  identified up to $ds\otimes d\P\otimes\lambda(da)$-a.e. equality, $K$ $\in$ $\calk^2$,   the subspace of $\cals^2$ consisting of nondecreasing processes  such that $K_0=0$, $\P$-a.s., and  the equation in  \eqref{BSDEnonpos} holds $\P$-a.s., while  the nonpositive jump constraint holds on $[0,T]\times\Omega\times A$ a.e. with respect to the measure $ds\otimes\P(d\omega)\otimes \lambda(da)$.

\begin{Definition}
A minimal solution to the BSDE \eqref{BSDEnonpos} is a quadruple $(Y,Z,U,K)$ $\in$ $\cals^2\times L^2(W)\times L^2(\tilde\mu)\times\calk^2$ solution to  \eqref{BSDEnonpos}  such that for any other solution
$(Y',Z',U',K')$ to the same BSDE, we have $\P$-a.s.:
\beqs
Y_s  & \leq & Y_s',  \;\;\; s \in [0,T].
\enqs
\end{Definition}

As noticed in Remark 2.1 in \cite{khph}, the minimal solution, when it exists, is unique as an element of $\cals^2\times L^2(W)\times L^2(\tilde\mu)\times\calk^2$
 By misuse of language,  we say sometimes that $Y$ (instead of the quadruple $(Y,Z,U,K)$) is the minimal solution  to \eqref{BSDEnonpos}.

 \vspace{3mm}

 We state the main result of this paper.

 \begin{Theorem} \label{main2}
For all  $(t,x,a)$ $\in$ $[0,T]\times\bfC^d\times A$, $\lambda$ $\in$ $\calm_f(A)$, and $\D_\lambda$ set-up, we have the following assertions:
\begin{enumerate}
\item There exists a unique minimal solution $Y$ $=$ $Y^{t,x,a,\lambda,\D_\lambda}$ to \eqref{BSDEnonpos}.  Moreover, for $s$ $\in$ $[0,t]$, $Y_s^{t,x,a,\lambda,\D_\lambda}$  is deterministic and does not depend on  $\D_\lambda$.
\item This minimal solution is related to the primal and dual control problems by:
\beq\label{indepalambda}
v(t,x) &=& v_1^*(t,x,a,\lambda,\D_\lambda)  \; = \;  v_2^*(t,x,a,\lambda) \; = \; Y_t^{t,x,a,\lambda,\D_\lambda}.
\enq
In particular, $v_1^*$ neither depends on $\D_\lambda$,  i.e. the strong and
weak dual control problems coincide, and none of the functions in \eqref{indepalambda} depends   on $(a,\lambda)$.
\end{enumerate}
\end{Theorem}

\begin{Remark}
{\rm  The HJB equation  for a stochastic control problem is usually derived from a dynamic programming principle on the value function.  We note that in view of the above relation
$v(t,x)$ $=$ $Y_t^{t,x,a,\lambda,\D_\lambda}$ (or simply
$v(t,x)$ $=$ $Y_t^{t,x}$ since it does not depend on $a,\lambda,\D_\lambda$), which is proved by purely probabilistic arguments, 
this gives another method to derive the HJB equation from the minimal BSDE solution $Y$. In the Markovian case, it is shown in \cite{khph} that the solution $Y^n$ to
 an appropriate penalized BSDE  satisfies a  semi linear PDE and, by passing to the limit, one obtains that $Y_t^{t,x}$  is solution to the nonlinear HJB equation. Such derivation does not resort to  dynamic programming principle, which is known to be a delicate issue, and would  be interesting to explore in the non-Markovian context.
\qed
}
\end{Remark}

\section{Proof of the main results}

\setcounter{equation}{0} \setcounter{Assumption}{0}
\setcounter{Theorem}{0} \setcounter{Proposition}{0}
\setcounter{Corollary}{0} \setcounter{Lemma}{0}
\setcounter{Definition}{0} \setcounter{Remark}{0}

\label{sectionproofmain}

Throughout this section, we make the standing assumptions {\bf (H1)} and {\bf (H2)}.

\subsection{Proof of Theorem \ref{main}}

We start with the inequality $v(t,x)$  $\ge$   $v_2^*(t,x,a,\lambda)$, which  is the easy part in the proof of Theorem \ref{main} and
it is to be expected since, intuitively,  in the dual control problem we control $X$ through the intensity of the control process $I$, so we cannot have a better performance than choosing
directly  the control process $\alpha$ as we do in the primal problem.

\subsubsection{Proof of the inequality $v$ $\ge$ $v_2^*$}

Fix $t$ $\in$ $[0,T]$, $x$ $\in$ $\bfC^d$, $a\in A$, $\lambda\in\calm_f(A)$ and an  admissible dual set-up   $\D_\lambda=(\Omega, \calf,  \P, W,\mu_\lambda)$.
Let   $\F$ $=$ $(\calf_t)_{t\ge 0}$ be  the completion of  the natural filtration generated by  $(W,\mu_\lambda)$.  Choose  a  random field $\nu$ in the class $\calv(\D_\lambda)$ and define the corresponding probability measure
$\P^\nu$. Define $(I,X)$, as the solution to \eqref{defxi} and consider
the gain for the dual control problem
\beqs
J^*(t,x,a,\lambda,\D_\lambda, \nu) &=& \E^\nu \Big[ \int_t^Tf_s(X,I_s)\,ds+g(X)\Big].
\enqs
Now we note that $\A$ $:=$ $(\Omega, \calf,\F, \P^\nu, W)$ is an admissible set-up for the original control problem. Moreover, since
the process $I$ is   progressive  with respect to  $\F$, it belongs to the class of $\A$-admissible controls, i.e. it is an element of $\cala(\A)$.
The corresponding trajectory can be obtained solving equation \eqref{stateeq} with $\alpha=I$, which coincides with the
equations defining the process $X$ in \eqref{defxi}. It follows that, with this choice of $\A$ and $\alpha$, the gain for the original and the dual problems are the same, so we have
 \beqs
 J^*(t,x,a,\lambda,\D_\lambda, \nu) &=& J(t,x,\A,I) \; \le \; v(t,x).
 \enqs
The required conclusion follows by taking the supremum  with respect to  $\nu \in \calv(\D_\lambda)$,  and then with respect to $\D_\lambda$,  in the left-hand side.

\subsubsection{Proof of the inequality  $v_2^*$ $\geq$ $v$}

The required inequality  will be a  consequence of the following proposition:

\begin{Proposition} \label{mainproposition}
Fix   $t\in [0,T]$, $x\in\bfC^d$, $a\in A$, $\lambda\in\calm_f(A)$.
Then, for every admissible set-up $\A'=(\Omega', \calf',\G', \Q', W')$, for   every admissible control
$\alpha\in\cala(\A')$ and  for every number $\delta>0$, there exist an  admissible dual set-up
$\D_\lambda=(\Omega, \calf,  \P, W,\mu_\lambda)$ and an element $\nu\in\calv(\D_\lambda)$
such that
\beqs
J^*(t,x,a,\lambda,\D_\lambda, \nu) &\ge&  J(t,x,\A',\alpha)-\delta.
\enqs
\end{Proposition}

We first check  that  the required inequality $v_2^*(t,x,a,\lambda)$ $\ge$  $v(t,x)$ is an  immediate consequence.
Indeed, from the inequality of Proposition \ref{mainproposition}, it follows that
 \beqs
 v_2^*(t,x,a,\lambda) &\ge& J(t,x,\A',\alpha) - \delta.
\enqs
Taking the supremum  with respect to  $\alpha \in \cala(\A')$,  and then with respect to all admissible set-ups $\A'$,  we conclude that
 $v_2^*(t,x,a,\lambda)$ $\ge$  $v(t,x)-\delta$ and finally the required inequality  follows from the arbitrariness of $\delta$.

 \vspace{3mm}

The rest of this section is devoted to the proof of Proposition \ref{mainproposition}.
Since the proof is rather technical, in order to simplify the notation we will limit the exposition to the case when $t=0$. The general case can be proved in the same way, with slight and
obvious changes. We fix  elements $x\in\bfC^d$, $a\in A$, $\lambda\in\calm_f(A)$.
We also fix an admissible set-up $\A'=(\Omega', \calf',\G', \Q', W')$,
an admissible control $\alpha\in\cala(\A')$, and $\delta>0$. The corresponding trajectory $X$ is the unique solution to
\begin{equation}\label{stateeqbis}
 X_t \; = \; x(0)+\int_0^t b_s(X,\alpha_s)\,ds  +\int_0^t      \sigma_s(X,\alpha_s)\,dW_s,\qquad t \in [0,T].
\end{equation}
Note that $x(\cdot)$ affects the trajectory $X$ only through the value $x(0)$.
The  gain functional is
\beqs
J(0,x,\A',\alpha) &=& \E^{\Q'}\, \Big[ \int_0^Tf_t(X,\alpha_t)\,dt+g(X)\Big].
\enqs

Following \cite{kr}, we introduce a metric in the set $\cala(\A')$ of admissible controls. Recall that $\rho$ denotes the metric
in $A$, chosen such that $\rho$ $<$ $1$.
%and note that $(x,y)\mapsto d(x,y)/(1+d(x,y))$, ($x,y\in A$) is an equivalent metric which is also bounded.
Next define, for any $\alpha^1,\alpha^2\in\cala(\A')$,
\beqs
\tilde \rho(\alpha^1,\alpha^2) &=& \E^{\Q'} \Big[\int_0^T\rho(\alpha^1_t,\alpha^2_t)\,dt \Big].
\enqs
Note that a sequence $\alpha^n$ converges in $\cala(\A')$ to a limit $\alpha$ with respect to this metric if and only if $\alpha^n\to\alpha$
in $dt\otimes d\Q'$-measure, i.e. if and only if
\beqs
\lim_{n\to\infty}(dt\otimes d\Q')(\{(t,\omega')\in[0,T]\times\Omega'\,:\, \rho(\alpha^n_t(\omega'),\alpha_t(\omega'))>\epsilon\}) &=& 0,
\qquad {\rm for\; any \;} \epsilon>0.
\enqs

In \cite{kr},  the following continuity result  of the gain functional with respect to the control  is proved in the case of controlled diffusion processes.
The extension to our non-Markovian situation is straightforward, so we only sketch its proof.

\begin{Lemma}\label{Jcontinuous}
The map $\alpha\mapsto J(0,x,\A',\alpha)$ is continuous with respect to the metric $\tilde \rho$.
\end{Lemma}
{\bf Proof.} In this proof we write $\E'$ instead of $\E^{\Q'}$ for short. Suppose $\alpha^n,\alpha\in\cala(\A')$ and  $\alpha^n\to\alpha$
in $dt\otimes d\Q'$-measure. Denote $X^n, X$ the corresponding
trajectories. Then, starting from the state equation \eqref{stateeqbis}, using usual arguments involving the Burkholder-Davis-Gundy inequalities
and the Gronwall lemma, for every $p\in [1,\infty)$ we arrive at
\beqs
\E' \Big[\sup_{t\in[0,T]}|X^n_t-X_t|^p \Big] &\le&  C \Big\{  \E' \Big[ \int_0^T | b_t(X,\alpha^n_t)-b_t(X,\alpha_t)|^p  + | \sigma_t(X,\alpha^n_t)-\sigma_t(X,\alpha_t)|^p  \,dt \Big] \Big\},
\enqs
for a suitable constant $C$, independent of $n$.  Recalling the bound \eqref{solusp} on the solution $X$,
by standard arguments we first conclude, by the dominated convergence theorem,
under {\bf (H1)}, that $\E'\Big[\sup_{t\in[0,T]}|X^n_t-X_t|^p\Big]$ $\to 0$ as $n\to\infty$.
Next we have
\beqs
|J(0,x,\A',\alpha^n)-J(0,x,\A',\alpha)| &\le&
\E' \Big[ \int_0^T|f_t(X^n,\alpha^n_t)- f_t(X,\alpha_t)|\,dt \Big] + \E'|g(X^n)-g(X)|.
\enqs
To finish the proof we show that the right-hand side tends to zero.
Suppose on the contrary,   that there exist $\eta>0$ and a subsequence (denoted $(X^{n'},\alpha^{n'})$) such that
\begin{equation} \label{contradictiondominated}
\E' \Big[\int_0^T|f_t(X^{n'},\alpha^{n'}_t)- f_t(X,\alpha_t)|\,dt \Big] \; \ge \;  \eta,
\end{equation}
for every $n'$. Passing to a sub-subsequence, still denoted by the same symbol,  we can assume that
\beqs
\sup_{t\in[0,T]}|X^{n'}_t-X_t|\to 0, \;\;\;d\Q'-a.s.,\quad
\rho(\alpha^n_t,\alpha_t)\to 0,\;\;\; dt\otimes d\Q'-a.e.
\enqs
as $n'\to\infty$, and by the assumed continuity properties of $f$ it follows that
$f_t(X^{n'},\alpha^{n'}_t)\to f_t(X,\alpha_t)$,  $dt\otimes d\Q'$-a.e.
Next we extract a further subsequence $(n'_j)$ such that
\beqs
\left(\E' \Big[ \sup_{t\in[0,T]}|X^{n'_j}_t-X_t|^p \Big] \right)^{1/p} &\le&  2^{-j},
\enqs
so that the random variable $\bar X:=\sum_j \sup_{t\in[0,T]}|X^{n'_j}_t-X_t|$ satisfies $\E'|\bar X|^p<\infty$ as well as
$|X^{n'_j}_t|\le |X_t|+ |\bar X|$ for every $t$ and $j$. Recalling the polynomial growth condition  of $f$ in {\bf (H2)}, we obtain
\beqs
|f_t(X^{n'_j},\alpha^{n'_j}_t)- f_t(X,\alpha_t)| & \le &  C(1+\sup_{t\in[0,T]}|X^{n'_j}_t|^m
+ \sup_{t\in[0,T]}|X_t|^m) \\
&\le &  C(1+ |\bar X|^m+ \sup_{t\in[0,T]}|X_t|^m)
\enqs
for a suitable constant $C$, and choosing $p$ $=$ $m$,  we conclude that the right-hand side is integrable, which gives a contradiction with
\eqref{contradictiondominated} by the dominated convergence theorem. This shows that
$\E' \big[\int_0^T|f_t(X^n,\alpha^n_t)- f_t(X,\alpha_t)|\,dt\big]$ $\to 0$, and in a similar way one shows that $\E'|g(X^n)- g(X)|\to 0$.
\qed

\begin{Remark}\begin{em}
\label{extendedcont}
For further use we note that the metric $\tilde \rho$ can be defined on the set of all $\calb([0,T])\otimes\calf'$-measurable
$A$-valued processes. Now suppose that we have a sequence of filtrations $\H^m$ $=$ $(\calh^m_t)$, satisfying $\calg'_t\subset  \calh^m_t$ for every $t\ge0$, such that $W$ is a Wiener process with respect to each of them, and denote $\A^m=(\Omega',\calf',\H^m, \Q',W)$ the corresponding set-ups. Then, given a sequence $\alpha^m\in\cala(\A^m)$, and $\alpha\in\cala(\A')$, the convergence $\tilde \rho(\alpha^m,\alpha)\to 0$ still implies $J(0,x,\A^m,\alpha^m)\to J(0,x,\A',\alpha)$. This is a slight extension of Lemma  \ref{Jcontinuous},
that can be proved by the same arguments as before.
\qed
\end{em}
\end{Remark}

\vspace{2mm}

The following lemma shows that we can replace any control process by
another control which is a pure jump process, without changing the gain too much.

\begin{Lemma}\label{reductiontopointcontrol}
There exists an admissible control  process $\bar\alpha\in\cala(\A')$,    such that
\begin{equation}\label{acsquasiottimobis}
   J(0,x,\A',\bar\alpha) \; \ge \;  J(0,x,\A',\alpha) -\delta,
\end{equation}
and $\bar\alpha$ has the form $\bar\alpha_t=\sum_{n= 0}^{N-1} \alpha_{n}1_{ [T_n,T_{n+1})}(t)$,
where  $0=T_0<T_1<\ldots T_N=T$ is a deterministic partition of $[0,T]$,
$\alpha_0,\ldots,\alpha_{N-1}$ are $A$-valued
random variables that take only a finite number of values, and each $\alpha_n$ is $\calg'_{T_n}$-measurable. In addition, we can choose $\bar\alpha$ satisfying $\alpha_0=a$.
\end{Lemma}
{\bf Proof.} In \cite{kr} Lemma 3.2.6 it is proved that the set of admissible
controls $\bar\alpha$ having the form specified in the lemma are dense in $\cala(\A')$ with respect to the metric $\tilde \rho$. The Lemma
is then a consequence of the continuity property of Lemma \ref{Jcontinuous}.
The additional requirement that $\alpha_0=a$ can be fulfilled by adding, if necessary, another point $T'$  close to $0$ to the subdivision
and modifying $\bar\alpha$ setting $\bar\alpha_t=a$  for   $t\in [0,T')$. This modification is as close as we wish to the original control with respect
to the metric $\tilde \rho$, provided $T'$ is chosen sufficiently small.
\qed

\vspace{2mm}

Before proceeding further, we need to recall some useful facts, summarized
in the following remark.

\begin{Remark}\begin{em}
\label{spaceenlargement}
We start with an admissible set-up $\A'=(\Omega', \calf',\G', \Q', W)$,  and we need to give new definitions  or make suitable constructions using
additional, \emph{independent}, random variables or stochastic processes.  If these random objects are not already
defined on $(\Omega', \calf',\Q')$ we can perform the following standard construction in order to reduce to this case. We consider
another probability space $(\Omega'', \calf'', \Q'')$ on which are defined these random objects, and we define   $(\Omega,\calf,\Q)$ setting
\beqs
\Omega \; = \; \Omega'\times \Omega'',
\qquad
\calf \; = \; \calf'\otimes \calf'',
\qquad
\Q \; = \; \Q'\otimes \Q''.
\enqs
We can also define a filtration $\G$ $=$ $(\calg_t)$ in $(\Omega,\calf)$ setting $\calg_t=\{B\times \Omega''\,:\, B\in \calg'_t\}$ for $t\ge 0$.

Any random variable $Z$ (respectively, stochastic process  $(Y_t)$)
  in $(\Omega',\calf')$
admits a natural extension to a   random variable (resp. stochastic process)
 on $(\Omega,\calf)$,
 still denoted by the same symbol, given by $Z(\omega)=Z(\omega')$
(resp. $Y_t(\omega)=Y_t(\omega')$), for  $\omega=(\omega',\omega'')\in\Omega$.
It easily verified that if $Y$ is $\G'$-adapted (resp. optional, progressive, predictable) then its extension is $\G$-adapted (resp. optional, progressive, predictable).
Moreover, the extension of $W$ is a Wiener process with respect to $\G$ and $\Q$, so that we have constructed another
admissible set-up  $\A:=(\Omega , \calf ,\G, \Q,W)$. Note that the extension of an
 $\A'$-admissible control is $\A$-admissible, that the trajectory of an extended control
 process is the extension of the original trajectory, and that the corresponding
 gain functional has the same value.

Similar considerations hold for random objects originally defined on $\Omega''$.
For instance, if $(V_n)$ denotes a random sequence on $\Omega''$,
its extension has the same law and is independent of $\calg_\infty$ (hence
independent of $W$) under $\Q$.

We shall briefly describe this construction by saying that $\Omega$ is an
{\em enlargement} of $\Omega'$, or that the set-up  $(\Omega , \calf ,\G, \Q,W)$ is an {\em enlargement} of $(\Omega' , \calf' ,\G', \Q',W)$.

In conclusion,  starting with our admissible set-up $\A'=(\Omega' , \calf' ,\G', \Q',W')$, we have proved the existence of an enlargement
 $\A=(\Omega , \calf ,\G, \Q,W)$ on which there exist   random elements  with arbitrary prescribed laws and independent of $\calg_\infty$ under $\Q$.
\qed
\end{em}
\end{Remark}

\vspace{2mm}

In the proof of Proposition \ref{mainproposition},  we need a preliminary result, stated below as Lemma \ref{compensatorofmodifiedppmbis}, where we basically prove
that the marked point process $\bar\alpha$ in the statement of Lemma \ref{reductiontopointcontrol} can be slightly perturbed in such a way that
its compensator becomes absolutely continuous with respect to the measure $\lambda(da)\,dt$. Then we shall see  in Lemma \ref{reductiontopointcontrolabscont}
that the corresponding gain will also be changed slightly. The proof depends on a more general and technical result, reported in the Appendix as Lemma \ref{MPPperturbed}.

Let us come back again to the original set-up $\A'=(\Omega' , \calf' ,\G', \Q',W')$ in the statement of Proposition \ref{mainproposition} and denote by
$\bar\alpha$ the admissible control of Lemma \ref{reductiontopointcontrol}.
Let  $\A=(\Omega , \calf , \G, \Q,W)$ be an  enlargement  of the set-up $\A'$,
as described in  Remark \ref{spaceenlargement}, and denote the extension of the control $\bar\alpha$ still by the same symbol.
It is convenient to extend further the definition of $\bar\alpha$ to $[0,\infty)\times\Omega$ in a trivial way setting
\beqs
\alpha_n &= & \alpha_{N-1}, \;\; n\ge N;\qquad T_n=T+n-N,
\;\; n>N,
\enqs
and defining $\bar\alpha_t=\sum_{n= 0}^{\infty} \alpha_{n}1_{ [T_n,T_{n+1})}(t)$.
This way, $\bar\alpha$ is associated to the marked point process $(T_n,\alpha_n)_{n\ge 1}$ and $\alpha_0=a$.

For every integer $m\ge1$, on $\Omega$  we can find sequences $(U_n^m)_{n\ge1}$, $(S^m_n)_{n\ge1}$ of real random variables satisfying the following conditions:
\begin{enumerate}
\item every $ U_n^m $ is uniformly  distributed on $(0,1)$;
\item  every   $ S_n^m$ admits
 a density (denoted $f_n^m(t)$)  with respect to the Lebesgue measure,
 and we have $ 0<S_1^m<S_2^m<S^m_3<\ldots $ and $\sum_{n\ge1}S_n^m<1/m$
 for every $m$;
 \item for every $m$, the random variables $U^m_n,S^m_k$ ($n,k\ge1$)  are independent, and independent  of $\calg_\infty$.
 \end{enumerate}

For every   $m\ge 1$, let $\bfB(b,1/m)$ denote the open ball  of radius $1/m$, with respect to the metric $\rho$, centered at $b\in A$.
Since $\lambda(da)$ has full support, we have $\lambda(\bfB(b,1/m))>0$ and we can  define a transition kernel $q^m(b,da)$ in $A$ setting
\beqs
q^m(b,da)&=& \frac{1}{\lambda(\bfB(b,1/m))}\, 1_{\bfB(b,1/m)}(a) \lambda(da).
\enqs
We recall that  we require $A$ to be a Lusin space. It follows from Lemma \ref{transkernel} that there exists a function $q^m:A\times [0,1]\to A$, measurable with respect to
$\calb(A)\otimes \calb([0,1])$,  such that for every $b\in A$ the measure $B\mapsto q^m(b,B)$ ($B\in \calb(A)$) is the image of the Lebesgue measure  on $[0,1]$ under the mapping
$u\mapsto q^m(b,u)$. Thus, if  $U$ is a random variable defined on some probability space and having uniform law on $[0,1]$ then, for fixed $b\in A$, the random variable
$q^m(b,U)$ has law $q^m(b,da)$ on $A$. The use of the same symbol $q^m$ should not generate confusion.

\vspace{1mm}

Define
 \begin{equation}\label{defRbeta}
    R^m_n \; = \; T_n+S^m_n,\quad \beta_n^m=q^m(\alpha_n,U^m_n),
    \qquad \qquad n,m\ge 1,
 \end{equation}
and set $R^m_0=0$. Since we assume $S^m_n<S^m_{n+1}$ we see that $(R^m_n,\beta^m_n)_{n\ge1}$
is a marked point process in $A$. Let $\kappa^m=\sum_{n\ge1}\delta_{(R^m_n,\beta^m_n)}$
denote the corresponding random measure, and $(\calf^{\kappa^m}_t)$ the corresponding natural filtration. Finally set
\begin{equation}\label{defaccam}
\calh_t^m \; = \; \calg_t\vee \calf^{\kappa^m}_t, \qquad t\ge 0.
\end{equation}
Now Lemma \ref{MPPperturbed} in the Appendix provides us with the following explicit formula for the compensator
$\tilde\kappa^m$ of $\kappa^m$ with respect to the filtration $\H^m$ $=$ $(\calh^m_t)$:
 \beqs
 \tilde\kappa^m(dt,da) &=& \sum_{n\ge 1}  1_{(T_n\vee R^m_{n-1},R^m_n]}(t)\,
 q^m(\alpha_n,da) \frac{f_n^m(t-T_n)}{1-F_n^m(t-T_n)}\,dt,
 \enqs
where we denote  by  $F_n^m(s)=\int_{-\infty}^sf_n^m(t)dt$ the cumulative distribution function of $S_n^m$, with  the convention that $\frac{f_n^m(s)}{1-F^m_n(s)}=0$ if $F^m_n(s)=1$.

\vspace{2mm}

We summarize the relevant properties of this construction in the following result.

\begin{Lemma}\label{compensatorofmodifiedppmbis}
With the previous notation, in the enlarged set-up $\A$ the following properties hold true:
\begin{enumerate}
  \item $T_n<R_n^m$ and $\sum_{n\ge1}(R_n^m-T_n)<1/m$;
  \item $\rho(\alpha_n,\beta^m_n)<1/m$;
  \item   the $(\Q,\H^m)$-compensator of   $\kappa^m$ is absolutely continuous with respect to $\lambda(da)\,dt$, so that it can be written in the form
  $$
 \tilde\kappa^m(dt,da) \; = \;  \phi^m_t(a)\, \lambda(da)\,dt
 $$
for a suitable nonnegative $\calp(\H^m)\otimes \calb(A)$-measurable  function $\phi^m$.
\end{enumerate}
%Here $\calp(\H^m)$ denotes the predictable $\sigma$-algebra with respect to the filtration $\H^m$ $=$ $(\calh_t^m)$.
\end{Lemma}
{\bf Proof.}
The first property follows from the fact that $S_n^m>0$ and  $\sum_{n\ge1}S_n^m<1/m$. Since, for every $b$ $\in$ $A$ , $q^m(b,da)$ is supported in $\bfB(b,1/m)$ we have
$\rho(\alpha_n,\beta^m_n)<1/m$.  Finally, from the choice of the kernel $q^m(b,da)$ we obtain
\beqs
\phi^m_t(a) &=& \sum_{n\ge 1}  1_{(T_n\vee R^m_{n-1},R^m_n]}(t)\,
\frac{1}{\lambda(\bfB(\alpha_n,1/m))}\, 1_{\bfB(\alpha_n,1/m)}(a)
  \frac{f_n^m(t-T_n)}{1-F_n^m(t-T_n)}.
\enqs
\qed

\vspace{2mm}

Now recall that we have fixed from the beginning $t=0$, $x\in\bfC^d$, $a\in A$, $\lambda\in\calm_f(A)$, a set-up $\A'$, $\alpha$ $\in$ $\cala(A')$ and $\delta$ $>$ $0$.
Also recall the notation $J(0,x,\A',\alpha)$ for the gain functional.

\begin{Lemma}\label{reductiontopointcontrolabscont}
%For every admissible set-up $\A'=(\Omega', \calf',\G', \Q', W')$, for   every admissible control $\alpha\in\cala(\A')$ and  for every $\delta>0$,
There exists an admissible set-up $\A''=(\Omega , \calf ,\H, \Q,W)$, which is  an enlargement  of $\A'$, and an admissible control $\hat\alpha\in \cala(\A'')$ such that
\begin{enumerate}
\item
   $J(0,x,\A'',\hat\alpha) > J(0,x,\A',\alpha)-2\delta;$
\item
there exists an $\H$-marked point process $(R_n,\beta_n)_{n\ge 1}$ such that
$
\hat\alpha_t=\sum_{n\ge 0} \beta_{n}1_{ [R_n,R_{n+1})}(t)
$ ($R_0=0,\beta_0=a$);
\item
  the $(\Q,\H)$-compensator of the corresponding random measure
  $\kappa=\sum_{n\ge1}\delta_{(R_n,\beta_n)}$
  is absolutely continuous with respect to $\lambda(da)\,dt$, so that
it can be written in the form
  $$
 \tilde\kappa(dt,da) \; = \;  \phi_t'(a)\, \lambda(da)\,dt
 $$
 for a suitable nonnegative $\calp(\H)\otimes \calb(A)$-measurable  function $\phi'$.
\end{enumerate}
\end{Lemma}
{\bf Proof.} We first construct an appropriate enlargement $\A=(\Omega, \calf,\G,\Q,W)$ of $\A'$.
Then we take the control $\bar\alpha$ of Lemma \ref{reductiontopointcontrol} and we extend it to
$[0,\infty)\times\Omega$  as described above, so that it is associated to the  $\G$-marked point process $(T_n,\alpha_n)_{n\ge 1}$.
 Finally, for every $m\ge 1$,  we introduce  $(R_n^m,\beta^m_n)_{n\ge1}$ and the filtration $\H^m$ $=$ $(\calh^m_t)$ defined by \eqref{defRbeta}
and \eqref{defaccam}.  Since the random variables $S^m_n$, $U^m_n$ occurring in \eqref{defRbeta} are independent of $\calg_\infty$
it follows that  $W$ is a Wiener process with respect to $\H^m$. Therefore  $\A^m:=(\Omega, \calf,\H^m,\Q,W)$ is an admissible
 set-up.  Next we define
 $$
\hat\alpha_t^m \; = \; \sum_{n\ge 0} \beta^m_{n}1_{ [R^m_n,R^m_{n+1})}(t)
$$
with the convention $R^m_0=0,\beta^m_0=a$, and note that it is an admissible control, i.e. an element of $\cala(\A^m)$.
Now let us compare those controls with $\bar\alpha_t=\sum_{n= 0}^{\infty} \alpha_{n}1_{ [T_n,T_{n+1})}(t)$.
The first two conclusions of Lemma \ref{compensatorofmodifiedppmbis} show that $\hat\alpha_t^m$ converges to $\bar\alpha$
in $d\Q\otimes dt$-measure as $m\to\infty$, hence with respect to the metric $\tilde \rho$ introduced before. By Lemma  \ref{Jcontinuous} and
Remark \ref{extendedcont}, this shows that $J(0,x,\A^m,\hat\alpha^m)\to J(0,x,\A',\bar\alpha)$ as $m\to\infty$.
So there exists $M$ so large such that
$J(0,x,\A^M,\hat\alpha^M) > J(0,x,\A',\bar\alpha)-\delta$ and,
by \eqref{acsquasiottimobis},
$J(0,x,\A^M,\hat\alpha^M) > J(0,x,\A',\bar\alpha)-2\delta$. We finally set
$$
\A'':=\A^M,  \quad \calh_t:=\calh^M_t,\quad
\beta_{n}:=\beta^M_{n},\quad  R_n:=R^M_n,
\quad \hat\alpha_t:=
\hat\alpha_t^M, \quad
\phi_t'(a):=\phi^M_t(a).
$$
\qed

\vspace{3mm}

For the rest of the proof of Proposition \ref{mainproposition}, only Lemma \ref{reductiontopointcontrolabscont}
and Lemma  \ref{Jcontinuous} will be used.  The idea  is now to add to the control an additional {\em independent} Poisson measure with
compensator $k^{-1}\lambda(da)\,dt$, with $k$ large, i.e. with  intensity so small that the gain is not much affected.
The formal construction is as follows.

Let $\A''$ $=$ $(\Omega , \calf ,\H, \Q,W)$, $(R_n,\beta_n)_{n\ge 0}$, $\hat\alpha_t$ $=$ $\sum_{n\ge 0} \beta_{n}1_{ [R_n,R_{n+1})}(t)$,
$\kappa$ $=$ $\sum_{n\ge1}\delta_{(R_n,\beta_n)}$ with compensator  $\tilde\kappa(dt,da)= \phi_t'(a)\, \lambda(da)\,dt$, denote the objects introduced in  Lemma \ref{reductiontopointcontrolabscont}.  By Remark \ref{spaceenlargement} we can assume that, for every integer $k\ge 1$, on the probability space $(\Omega, \calf, \Q)$ there exists a Poisson random measure on $(0,\infty)\times A$, denoted
\beqs
\pi^k &=& \sum_{n\ge 1}\delta_{(T_n^k , {\xi_n^k })},
\enqs
admitting compensator $k^{-1}\lambda(da)\,dt$ with respect to its natural filtration $\F^{\pi^k}$ $=$ $(\calf^{\pi^k}_t)$,
and independent
of $\calh_\infty$.
Now we define another random measure setting
$$
\mu^k =\kappa+\pi^k .
$$
Note that the jumps times $(R_n)$  are independent
of the jump times $(T_n^k )$, and the latter have
absolutely continuous laws. It follows that, except possibly
on a set of $\Q$ probability zero, their graphs  are
disjoint, i.e. $\kappa$ and $\pi^k $ have no common jumps. Therefore,
  the random measure $\mu^k $ admits a representation
$$
\mu^k =\sum_{n\ge 1}\delta_{(S_n^k , {\eta^k _n})}
$$
where $(S_n^k ,\eta^k _n)_{n\ge 1}$ is a marked point process,  each $S_n^k $ coincides with one of the times $R_n$ or one of the times
$T_n^k $, and each $\eta_n^k $ coincides with one of the random variables $\xi_n^k $ or one of the random variables $\beta_{n}$.
Since $\kappa$ and $\pi^k $ are independent it follows    from Proposition \ref{sumofindepmeasures}   that
 $\mu^k $ has compensator  $(\phi'_t(a)+k^{-1})\,\lambda(da)\,dt$ with respect
 to the filtration $\H\vee\F^{\pi^k}$ $=$ $(\calh_t\vee\calf^{\pi^k }_t)$.
Let us  denote by $\F^k$ $=$ $(\calf_t^k )$  the completion of the natural filtration generated by $(W,\mu^k )$.
Clearly,  $W$ is a Wiener process with respect to $\F^k$.

Now we need to prove that the compensator of $\mu^k$ with respect to $\F^k$ remains absolutely continuous with respect to
$\lambda(da)dt$, see  Lemma \ref{compensabscont} below. For its formulation, consider the measure space
$([0,\infty)\times\Omega\times A, \calb([0,\infty))\otimes\calf\otimes \calb(A),dt\otimes\Q(d\omega)\otimes\lambda(da))$.
Although this is not a probability space, one can define in a standard way the conditional expectation
 of any positive measurable
function, given an arbitrary sub-$\sigma$-algebra.
Let us denote by  $\phi_t(\omega,a)$ the conditional expectation  of the  random field $ \phi_t'(\omega,a)$
 with respect to the $\sigma$-algebra $\calp(\F^k)\otimes \calb(A)$.

\begin{Lemma}\label{compensabscont}
The compensator of $\mu^k $ with respect to $(\Q,\F^k)$ is $(\phi_t(a)+k^{-1})\,\lambda(da)\,dt$.
\end{Lemma}
{\bf Proof.} Let $H$ denote an arbitrary positive $\calp(\F^k)\otimes \calb(A)$-measurable function.
Denote by $\F^{0,k}$ $=$ $(\calf_t^{0,k} )$  the
uncompleted  natural filtration generated by $(W,\mu^k )$.
Then there exists a
positive $\calp(\F^{0,k})\otimes \calb(A)$-measurable function $H'$ such that,
for $\Q$-almost all $\omega$, $H_t(\omega,a)=H'_t(\omega,a)$ for every $t,a$.
Since clearly
$\calf_t^{0,k} \subset \calh_t\vee\calf^{\pi^k }_t$ for every $t\ge0$,
 $H'$ is also measurable with respect to
$\calp(\H\vee\F^{\pi^k})\otimes \calb(A)$ and we have
\begin{eqnarray*}
\E^\Q \Big[ \int_0^\infty\int_AH_t(a)\,\mu^k(dt,da) \Big]
&=& \E^\Q \Big[ \int_0^\infty\int_AH_t'(a)\,\mu^k(dt,da) \Big]\\
&=& \E^\Q \Big[ \int_0^\infty\int_AH_t'(a)\,( \phi'_t(a)+k^{-1})\,\lambda(da)\,dt \Big]  \\
&=& \E^\Q \Big[ \int_0^\infty\int_AH_t(a)\,( \phi'_t(a)+k^{-1})\,\lambda(da)\,dt\Big]\\
&=& \E^\Q \Big[ \int_0^\infty\int_AH_t(a)\,( \phi_t(a)+k^{-1})\,\lambda(da)\,dt\Big],
\end{eqnarray*}
which proves the required result.
\qed

\vspace{3mm}

 Let us define the pure jump  process valued in  $A$  associated to the random measure $\mu^k$ by setting
 \beqs
 I^k _t &=& \sum_{n\ge 0} \eta^k _{n}1_{ [S_n^k ,S^k _{n+1})}(t), \;\;\; t \geq 0,
 \enqs
with the convention that $\eta^k _{0}=a$ ($a\in A$ was arbitrary and fixed above).
Then $\A^{k }$ $:=$ $(\Omega,\calf, \Q, $ $\F^k, W )$ is an admissible set-up and $I^k$ is an admissible control, i.e. $I^k\in\cala(\A^k)$.
We can prove that, for large $k$, $I^k$ is close to $\hat\alpha$ with respect to the metric $\tilde \rho$:

\begin{Lemma}\label{approxvalue} We have $\tilde \rho(I^k,\hat\alpha)\to 0$
as  $k\to\infty$.
\end{Lemma}
{\bf Proof.} We have  to prove that $I^k\to \hat\alpha$ in $dt\otimes d\Q$-measure.
Recall that the jump times of $\pi^k $ are denoted $T_n^k $.
Since  $T_1^k $ has exponential law with parameter
$ \lambda(A)/k$ the event $B_k =\{T_1^k  >T\}$ has probability
$e^{-\lambda(A)T/k}$, so that $\Q(B_k)\to 1$ as $k\to \infty$.
Noting that  on the set $B_k$, we have $\hat\alpha_t =I^k _t$ for all $t\in [0,T]$,
 the conclusion follows immediately. We remark that we have used the  fact that $\alpha_0=\eta_0^k=a$.
\qed

\vspace{3mm}

Applying Lemma \ref{Jcontinuous} and Remark \ref{extendedcont} we conclude
 that $J(0,x,\A^k,I^k)\to J(0,x,\A'',\hat\alpha)$ as $k\to\infty$.
So there exists $k$ large enough  such that $J(0,x,\A^k,I^k) > J(0,x,\A'',\hat\alpha)-\delta$ and,
by Lemma \ref{reductiontopointcontrolabscont},  $J(0,x,\A^k,I^k) > J(0,x,\A',\alpha)-3\delta$.
Introducing the notation
$$\A''':=\A^k,  \quad \calf_t:=\calf^k_t,\quad
\eta_{n}:=\eta^k_{n},\quad  S_n:=S^k_n,
\quad I_t:= I_t^k, \;\;\; t \geq 0,
$$
we see that we have proved the following result, where we choose $\epsilon=k^{-1}$
(recall that we have fixed from the beginning $t=0$, $x\in\bfC^d$, $a\in A$, $\lambda\in\calm_f(A)$, set-up $\A'$, $\alpha$ $\in$ $\cala(\A')$, and $\delta$ $>$ $0$):

\begin{Lemma}\label{reductiontopointcontrolabscontwitheps}
%For every admissible set-up $\A'=(\Omega', \calf',\G', \Q', W')$, for every admissible control $\alpha\in\cala(\A')$ and  for every $\delta>0$,
There exists an admissible set-up $\A'''=(\Omega , \calf ,\F, \Q,W)$, which is  an
 enlargement  of $\A'$, and an admissible control $I\in \cala(\A''')$ such that
\begin{enumerate}
\item
   $J(0,x,\A''',I) > J(0,x,\A',\alpha)-3\delta;$
\item
there exists an $\F$-marked point process $(S_n,\eta_n)_{n\ge 1}$
such that
$I_t=\sum_{n\ge 0} \eta_{n}1_{ [S_n,S_{n+1})}(t) $ ($S_0=0,\eta_0=a$);
\item $\F$ is the completion of the natural filtration of $(W,\mu)$, where
$\mu=\sum_{n\ge1}\delta_{(S_n,\eta_n)}$ is the corresponding random measure;
\item
 the $(\Q,\F)$-compensator of    $\mu$ is absolutely continuous 
 with respect to $\lambda(da)\,dt$ and
it can be written in the form
\beqs
 \tilde\mu(dt,da) &=&  (\phi_t(a)+\epsilon)\, \lambda(da)\,dt
\enqs
 for  some  $\epsilon>0$ and a nonnegative $\calp(\F)\otimes \calb(A)$-measurable  function $\phi$.
\end{enumerate}
\end{Lemma}

\vspace{3mm}

We note for further use that  the process $I$ and the corresponding trajectory $(X_t)_{t\in [0,T]}$ are the solution to
\begin{equation}\label{xicanonici}
\left\{
    \begin{array}{llll}
I_t &=&\dis a+\int_{(0,t]} \int_A (a'-I_{s-})\,\mu(ds,da'),
& t\in [0,T],
\\
 X_t&=&\dis x(0)+\int_0^t b_s(X,I_s)\,ds  +\int_0^t
      \sigma_s(X,I_s)\,dW_s,& t\in [0,T].
\end{array}
\right.
\end{equation}
which coincides with \eqref{defxi} in the case $t=0$ that we are addressing.

The final step in the proof of Proposition \ref{mainproposition} consists in showing that
the addition of the noise $\pi^k$  above
(a noise with intensity ``of size $\epsilon=k^{-1}$'') now makes it possible
to make a Girsanov transformation and construct
a dual admissible set-up where $\mu$ is a
Poisson random measure with compensator $\lambda(da)\,dt$
as required to fit the framework for the dual control problem described
in paragraph \ref{paradual}.

\vspace{3mm}

\noindent {\bf End of the proof of Proposition \ref{mainproposition}.}
Recall that we fix  $x\in\bfC^d$, $a\in A$, $\lambda\in\calm_f(A)$ and,
without loss of generality, $t=0$. We take arbitrary admissible set-up $\A'=(\Omega', \calf',\G', \Q', W')$, admissible control $\alpha\in\cala(\A')$ and  $\delta>0$.
Next we consider again the set-up  $\A'''=(\Omega , \calf ,\F, \Q,W)$,  the marked point process $(S_n,\eta_n)_{n\ge 1}$ (with $S_0=0, \eta_0=a$),
the corresponding admissible control $I_t$ $=$ $\sum_{n\ge 0} \eta_{n}1_{ [S_n,S_{n+1})}(t)$
and random measure   $\mu$ $=$ $\sum_{n\ge1}\delta_{(S_n,\eta_n)}$ and its compensator
$\tilde\mu(dt,da)$ $=$ $(\phi_t(a)+\epsilon)\, \lambda(da)\,dt$  in the statement of Lemma \ref{reductiontopointcontrolabscontwitheps},
and we recall that we have
\begin{equation}\label{perfinireottimo}
J(0,x,\A''',I) > J(0,x,\A',\alpha)-3\delta.
\end{equation}
We want to show that there exist
\begin{enumerate}
\item a probability measure $\P$ on $(\Omega,\calf)$ such that
\begin{equation}\label{defdualsetup}
    \D_\lambda=(\Omega, \calf,  \P, W,\mu)
\end{equation}
is  an   admissible dual set-up;
\item  an element $\bar\nu\in\calv(\D_\lambda)$ such that
\begin{equation}\label{dualequaldirect}
    J^*(0,x,a,\lambda,\D_\lambda, \bar\nu) \; > \; J(0,x,\A''',I)- \delta.
\end{equation}
\end{enumerate}
By \eqref{perfinireottimo}, and since $\delta$ is arbitrary,  this is enough to finish the proof  of Proposition \ref{mainproposition}.

\vspace{3mm}

To this end, let us define $\nu$ $=$ $\phi+\epsilon$,  and note  that $\nu$ is a positive $\calp(\F)\otimes\calb(A)$-measurable
random field. Since $\nu_t(a)\lambda(da)dt$ is the compensator of the non-explosive process $\mu$, it follows
easily that $\int_0^T\nu_t(a)\lambda(da)dt<\infty$ $\Q$-a.s., so that we can and will assume that $\nu$ takes only finite values on $[0,T]\times\Omega\times A$. Finally, since  the inverse $\nu^{-1}$ is  bounded and
$\lambda(A)$ $<$ $\infty$, it follows from standard criterion (see e.g. \cite{proshi08}) that the Dol\'eans-Dade exponential process
\begin{equation}\label{defgirsanovdensity}
 M_t \; := \; \exp\Big(\int_0^t\int_A(1-\nu_s(a)^{-1})\, \nu_t(a) \lambda(da)\,ds\Big) \prod_{S_n\le t} \nu_{S_n}(\eta_n)^{-1},\qquad t\ge 0,
\end{equation}
is a strictly positive martingale on $[0,T]$ (with respect to $\F$ and $\Q$), and we can
define an equivalent probability $\P$ on $(\Omega,\calf)$ setting
 $\P(d\omega)=M_T(\omega)\Q(d\omega)$. The expectation under
 $\P$ will be denoted simply $\E$ (while the expectation under  $\Q$ is denoted   $\E^\Q$).
Now we make the following claims:
\begin{enumerate}
  \item [(i)]$\mu$ (or more precisely its   restriction  to $(0,T]\times A$)
  has compensator $\lambda(da)\,dt$ (with respect to $\F$ and $\Q$);
  in particular,   it is a Poisson random measure.

  This follows from a theorem of  Girsanov type   (see \cite{ja} Theorem (4.5)),
  which guarantees that
  under the new probability $\P$   the   compensator of $\mu$   is given by
  $
  \nu_t(a)^{-1} \nu_t(a) \,\lambda(da)\,dt $ $=$ $ \lambda(da)\,dt.
  $
   \item [(ii)] $W$ is a $(\P,\F)$-Wiener process.

The proof is as follows. Since  the probabilities
$\P$ and $\Q$ are equivalent,  the quadratic variation of $X$ computed under $\P$ and $\Q$
is the same, and equals $\<W\>_t=t$. So it is enough to show
that $W$ is a $(\P,\F)$-local martingale, which is equivalent
to the fact that $MW$ is a $(\Q,\F)$-local martingale. Finally, this follows from a general fact: since $M$ is a $(\Q,\F)$-martingale
of finite variation, it is purely discontinuous  and
therefore   orthogonal (under $\Q$)
to $W$; thus, their product $MW$ is a $(\Q,\F)$-local martingale.

  \item  [(iii)]$W$ and $\mu$ are independent under $\P$.

To prove this claim it is enough to show that, for any measurable $B\subset A$, the process
 \beqs
  N^B_t &:=& \int_0^t\int_B \mu(ds,da) \; = \; \sum_n 1_{\{S_n\le t\}}1_{\{\eta_n\in B\}},
\enqs
is  independent from
  $W$  under $\P$.  From  claims (i) and (ii) it  follows that
  $N^B$ is a Poisson process and $W$ is a Wiener process, both with
  respect to $\F$ and $\P$. By a general result,
  see e.g. Theorem 11.43 in \cite{hewayan},  to check the independence it
  is enough to note that their right bracket $[N^B,W]$ is null,
  which is obvious, since $W$ is continuous and $N^B$ has
  no continuous part.
\end{enumerate}

From claims (i), (ii), (iii), and recalling that $\F$ is the completion of the natural filtration of $(W,\mu)$,  we deduce that
$\D_\lambda$ defined in \eqref{defdualsetup} is indeed an admissible  dual set-up. Note that we have checked that the $\P$-compensator of $\mu$
is $\lambda(da)dt$, although we do not make it explicit in the notation.

Next we proceed to verify \eqref{dualequaldirect}. Since  in general we can not assert that the random field
$\nu$ is bounded on $[0,T]\times\Omega\times A$, we can not conclude that
  it belongs to the class $\calv(\D_\lambda)$. However, we can still define
the process  $L^\nu$ by formula \eqref{LGirsanov}, which defines a strictly positive local martingale hence supermartingale with respect to $\P$.
It follows immediately from formulae \eqref{LGirsanov} and \eqref{defgirsanovdensity}  that
  $L_T^\nu$
is the inverse of $M_T$. It follows that
$\E [L_T^\nu]=\E_\Q [M_TL_T^\nu]=1$, so that
$L^\nu$ is indeed
a  $\P$-martingale on $[0,T]$ and   we can
define the corresponding probability $\P^\nu(d\omega):= L^\nu_T(\omega)\P(d\omega)$. Thus,
 the  Girsanov transformation $\P\mapsto \P^\nu$ is the inverse to the transformation $\Q\mapsto \P$ made above,
and changes back the probability $\P$ into $\P^\nu=\Q$ considered above.
In addition, we recall that the control $I\in \cala(\A''')$ constructed in Lemma \ref{reductiontopointcontrolabscontwitheps}
and the corresponding trajectory $(X_t)_{t\in [0,T]}$
are the solution to the system \eqref{xicanonici}, which
coincides with \eqref{defxi},
since we are assuming $t=0$. It follows that
\begin{equation}\label{gaindirectdual}
    J(0,x,\A''',I) \; = \;
\E^\Q \Big[ \int_0^Tf_t(X,I_t)\,dt+g(X)\Big]
\; = \; \E^\nu \Big[ \int_0^Tf_t(X,I_t)\,dt+g(X)\Big],
\end{equation}
where $\E^\nu$ denotes the expectation under $\P^\nu$.  If  $\nu$ belongs to $\calv(\D_\lambda)$, the right-hand side
equals the gain $J^*(0,x,a,\lambda,\D_\lambda, \nu)$  for the dual control problem and the desired inequality \eqref{dualequaldirect} obviously holds
with $\bar\nu=\nu$. However, since in general we can not assert that $\nu\in\calv(\D_\lambda)$,  we revert to the following approximation procedure by  truncation.

For any integer $k\ge 1$ define $\nu^k_t(a)=\nu_t(a)\wedge k$. Therefore
$\nu^k\in\calv(\D_\lambda)$,  we can define the corresponding process  $L^{\nu^k}$ by formula \eqref{LGirsanov},
 the probability  $\P^{\nu^k}(d\omega)=L_T^{\nu^k}(\omega)\,\P(d\omega)$, and
compute the gain
\beqs
J^*(0,x,a,\lambda,\D_\lambda, \nu^k) &=& \E^{\nu^k} \Big[ \int_0^Tf_t(X,I_t)\,dt+g(X)\Big],
\enqs
where $\E^{\nu^k} $ denotes the expectation under $\P^{\nu^k}$.
We claim that, for $k\to\infty$,
\begin{equation}\label{convgainsapprox}
    \E^{\nu^k} \Big[ \int_0^Tf_t(X,I_t)\,dt+g(X)\Big] \; \longrightarrow \;
\E^{\nu} \Big[ \int_0^Tf_t(X,I_t)\,dt+g(X)\Big].
\end{equation}
If we can prove the claim, it follows from
\eqref{gaindirectdual} that the inequality
\eqref{dualequaldirect} is verified with $\bar\nu=\nu^k$
and $k$ sufficiently large. So it remains to prove \eqref{convgainsapprox}, that
we re-write in the form
\begin{equation}\label{convgainsapproxbis}
\E^{\nu^k} [\Phi] \;=\;    \E [L^{\nu^k}_T\Phi] \;  \longrightarrow \;
\E \big[L^{\nu}_T \Phi\big] \;=\;  \E^{\nu} \big[ \Phi\big],
\end{equation}
where we have set $\Phi:= \int_0^Tf_t(X,I_t)\,dt+g(X)$. We note for further use that the assumption
of polynomial growth in {\bf (H2)}  implies that      there exists a constant $C$ such that
\begin{equation}\label{stimasuphipoly}
    |\Phi| \; \le \;    C(1+\sup_{t\in [0,T]}|X_t|^m).
\end{equation}
For $N\ge1$ we define $A_N=\{\sup_{t\in [0,T]}|X_t|>N\}$ and obtain
$$
\big|\E^{\nu^k} [\Phi]-\E^{\nu} [ \Phi]\big| \; \le \;
    \E^{\nu^k} [1_{A_N}|\Phi|]
    + \E^{\nu} [1_{A_N}|\Phi|] + \E [|L_T^{\nu^k}-L^{\nu}_T|1_{A_N^c}|\Phi|].
$$
By \eqref{stimasuphipoly} we have
$$
    \E^{\nu^k} [1_{A_N}|\Phi|] \; \le \;  C  \,   \E^{\nu^k} \Big[1_{A_N} (1+\sup_{t\in [0,T]}|X_t|^m)\Big] \; \le \;
\frac{C}{  N}\,   \E^{\nu^k} \Big[ \sup_{t\in [0,T]}|X_t| (1+\sup_{t\in [0,T]}|X_t|^m)\Big] \; \le \;
\frac{C'}{  N}
$$
for a suitable constant $C'$, independent of $k$, where for the last
inequality we have used the estimate
\eqref{solusp} whose right-hand side is the same for all
probabilities $\P^{\nu^k}$.
A similar estimate holds for $\E^{\nu} [1_{A_N}|\Phi|]$ and we obtain
\begin{equation}\label{stimaquasifinale}
    |\E^{\nu^k} [\Phi]-\E^{\nu} \left[ \Phi\right]|\le
\frac{2C'}{  N}+
\E [|L_T^{\nu^k}-L^{\nu}_T|1_{A_N^c}|\Phi|].
\end{equation}
By the dominated convergence theorem  we have
 $$
 \int_0^T\int_A(1-\nu_s^k(a))\, \lambda(da)\,ds
\; \longrightarrow\;  \int_0^T\int_A(1-\nu_s(a))\, \lambda(da)\,ds,
 $$
 a.s. (with respect to any of the equivalent
 probabilities $\P,\P^\nu,\P^{\nu^k}$)
 and
$$
 \prod_{S_n\le T} \nu^k_{S_n}(\eta_n) \; \longrightarrow\;
 \prod_{S_n\le T} \nu_{S_n}(\eta_n),
 $$
 a.s., since the product has finitely many factors a.s.  From formula
 \eqref{LGirsanov} that defines $L^\nu$ and $L^{\nu^k}$ we obtain
$L_T^{\nu^k}\to L_T^\nu$ a.s.
Since  $\E[L_T^{\nu^k}]=\E[L_T^\nu]=1$,
 we even have $L_T^{\nu^k}\to L_T^\nu$ in $L^1(\Omega,\calf,\P)$.
By \eqref{stimasuphipoly} the random variable
$1_{A_N^c}|\Phi|$ is bounded a.s.,
so letting $k\to\infty$ in
\eqref{stimaquasifinale} we have
$$
\limsup_{k\to\infty}  |\E^{\nu^k} [\Phi]-\E^{\nu} \left[ \Phi\right]|\le
\frac{2C'}{  N}.
$$
Letting $N$ tend to infinity  we  conclude  the proof of the claim \eqref{convgainsapproxbis},
and the proof of
Proposition \ref{mainproposition} is also finished.
\qed

\subsection{Proof of Theorem \ref{main2}}

We fix  initial conditions $(t,x,a)$ $\in$ $[0,T]\times\bfC^d\times A$,  $\lambda\in\calm_f(A)$, a dual set-up $\D_\lambda=(\Omega, \calf,  \P, W,\mu)$, denote by
$\F$ the completion of the natural filtration of $(W,\mu)$, and  consider the pair of processes $(I,X)$  as the solution to the system  \eqref{defxi}.  Note that we write $\mu$ instead of $\mu_\lambda$ for simplicity, and because  our results do not depend a posteriori on the choice of
$\lambda$.
Let us  introduce the family of penalized BSDEs associated to \eqref{BSDEnonpos}, parametrized by integer $n$ $\geq$ $1$:
\beq\label{bsdepenalized}
Y_s^n &=&  g(X) + \int_s^T f_r(X,I_r) dr  + n \int_s^T \int_A U_r^n(a)^+ \lambda(da) dr \\
 & &  -  \int_s^T Z^n_r\,dW_r  - \int_s^T\int_AU^n_r(a)\,\mu (dr,da), \;\;\; 0 \leq s \leq T,  \nonumber
\enq
where $u^+$ $=$ $\max(u,0)$.  It follows from a result in \cite{lita}, Lemma 2.4, that there exists a unique solution
$(Y^n,Z^n,U^n)\in \cals^2\times L^2(W)\times L^2(\tilde\mu)$ to \eqref{bsdepenalized}, where we drop the dependence on $(t,x,a,\lambda,\D_\lambda)$.
In that paper this result is proved by a cla\-ssical argument in the theory of
BSDE, namely a Picard iteration technique combined with a martingale representation theorem. This proof makes clear the following facts:

\begin{enumerate}
\item  The process $Y_s^n$ is deterministic for $s\in [0,t]$. Indeed,
denote by $\F^t=(\calf^t_s)_{s\ge t}$ the completed filtration on $[t,\infty)$
generated by the  restriction of the random
measure $\mu$ to $(t,\infty)\times A$ and
by the increments
of the Wiener process $W$ on $[t,\infty)$. Then, recalling that
  $X_s$ $=$ $x(s)$ and $I_s=a$ for $s\in [0,t]$, one sees that
$(X,I)$ is $\F^t$-progressive and 
that the BSDE \eqref{bsdepenalized} can be solved
on the   time interval $[t,T]$ with respect to $\F^t$. In particular it follows that $Y_t^n$ is $\calf^t_t$-measurable, hence deterministic. 
Setting, for $s\in [0,t]$,
$$
Z_s^n=0,\quad
U_s^n(a)=0,\quad
Y^n_s  =  Y^n_t +  \int_s^t f_r(x,a)\,dr, 
$$
one obtains the solution to \eqref{bsdepenalized} on $[0,T]$,
and $Y_s^n$ ($s\in [0,T]$)  is deterministic as claimed.
\item The law of the solution $(Y^n,Z^n,U^n)$, hence in particular the (deterministic) values of $Y_s^n$, for $s\in [0,T]$, is determined
by the coefficients $b,\sigma,f,g$ as well as the fixed elements $T, \lambda,a,t, x$ and $n$, but it does not depend on the particular choice of the probability space,
the Wiener process and the Poisson process. Thus, $Y^n_s$
($s\in [0,t]$) has the same value if equation \eqref{bsdepenalized} is solved
in another dual admissible control system $\D'_\lambda=(\Omega', \calf',  \P', W',\mu')$, provided $\mu'$ has $\P'$-compensator $\lambda(da)dt$ with respect to the natural
filtration generated by $(W',\mu')$.
\end{enumerate}

We provide an explicit representation of the solution to the penalized BSDE in terms of a family of auxiliary dual control problems. For every integer $n\ge1$,
let $\calv^n(\D_\lambda)$ denote the subset of elements $\nu_t(\omega,a)$ in $\calv(\D_\lambda)$ valued in $(0,n]$.  The following result is a  slight
modification of \cite{khph}, Proposition 4.1.

\begin{Lemma}\label{bsdepenlizaedandcontrol}
We have for all $n$ $\geq$ $1$,
\beq \label{expliYn}
Y_s^n &=& \mathop{\rm ess\, sup}_{\nu\in\calv^n(\D_\lambda)}  \E^\nu \Big[ \int_s^T f_r(X,I_r) \,dr +  g(X) \big| \calf_s \Big], \;\;\; s \in [0,T],  \qquad \P-a.s.
\enq
\end{Lemma}
{\bf Proof.} Fix $n$ $\geq$ $1$, and consider $(Y^n,Z^n,U^n)$ the solution to \eqref{bsdepenalized}.
For any $\nu\in\calv^n(\D_\lambda)$, let  $\mu^\nu(ds,da)$ $:=$ $\mu(ds,da)-\nu_s(a)\lambda(da)ds$ denote the $\P^\nu$-compensated martingale measure of $\mu$.
It is shown in \cite{khph}, Lemma 4.2 that for $Z^n$ $\in$ $L^2(W)$ and $U^n$ $\in$ $L^2(\tilde\mu)$, the processes
\beqs
\int_0^\cdot Z^n dW, & & \int_0^\cdot \int_A U_s^n(a) \mu^\nu(ds,da),
\enqs
are $\P^\nu$-martingales. Therefore, by taking the conditional expectation $\E^\nu$ given $\calf_s$ in  \eqref{bsdepenalized}, we obtain:
\begin{eqnarray*}
 Y_s^n &=& \E^\nu \Big[ \int_s^T f_r(X,I_r)\,dr + g(X) \big| \calf_s\Big]  \\
 & & + \; \E^\nu\Big[ \int_s^T\int_A [n  U^n_r(a)^ + -  \nu_r(a) U^n_r(a) ] \lambda(da) dr   \big| \calf_s \Big], \;\;\; s \in [0,T].
\end{eqnarray*}
From the elementary numerical inequality: $nu^+-\nu u$ $\geq$ $0$ for all $u$ $\in$ $\R$, $\nu$ $\in$ $[0,n]$, we deduce that
\begin{eqnarray*}
 Y_s^n & \geq & \mathop{\rm ess\, sup}_{\nu\in\calv^n(\D_\lambda)}  \E^\nu \Big[ \int_s^T f_r(X,I_r)\,dr + g(X) \big| \calf_s\Big],   \;\;\; s \in [0,T].
\end{eqnarray*}
For  $\epsilon$ $\in$ $(0,1)$,  define $\nu^\epsilon_s(a)$ $=$ $n1_{\{U^n_s(a)\ge 0\}}$ $+$ $\epsilon 1_{\{-1<U^n_s(a)< 0\}}$ $-$ $\epsilon U^n_s(a)^{-1} 1_{\{ U^n_s(a)\le-1\}}$.
Then  $\nu^\epsilon$ $\in$  $\calv^n(\D_\lambda)$,   and we have
\beqs
n U^n_s(a)^+  -  \nu_s^\epsilon(a) U^n_s(a) &\le & \epsilon, \;\;\;  0 \leq s \leq T,
\enqs
so that
\begin{eqnarray*}
Y_s^n &\le& \E^{\nu^\epsilon} \Big[ \int_s^T f_r(X,I_r)\,dr + g(X) \big| \calf_s\Big] + \epsilon T\lambda(A)  \\
 &\le& \mathop{\rm ess\, sup}_{\nu\in\calv^n(\D_\lambda)}  \E^\nu \Big[ \int_s^Tf_r(X,I_r) \,dr + g(X) \big| \calf_s\Big]+\epsilon T\lambda(A),
\end{eqnarray*}
which is enough to conclude the proof. Note that we could not take $\nu_s(a)$ $=$ $n 1_{\{U^n_t(a)\ge 0\}}$, since this process does not belong to $\calv^n(\D_\lambda)$
because of the requirement of strict positivity.
\qed

\vspace{3mm}

As a consequence of this explicit representation of  the penalized BSDE, we  obtain the follo\-wing uniform estimate on the sequence $(Y^n)_n$:

\begin{Lemma}\label{estimatesonbsdepenlizaed}
The sequence $(Y^n)_n$ is monotonically increasing in $n$, and we have
\beqs
\sup_{s\in [0,T]}|Y_s^n| &\le & C  \Big(1+\sup_{s\in [0,T]}|X_s|^m\Big),  \qquad \P-a.s.
\enqs
 for some constant $C$ depending only on  $T$, $m$ and on the constants $K_1,K_2$ as defined in {\bf (H1)}, {\bf (H2)}.
\end{Lemma}
{\bf Proof.} Monotonicity follows from the formula for $Y^n$ presented in Lemma  \ref{bsdepenlizaedandcontrol}, since $\calv^n(\D_\lambda)$ $\subset$ $\calv^{n+1}(\D_\lambda)$.
Then the inequality $Y^n_s$ $\le$ $Y^{n+1}_s$ holds $\P$-a.s. for all $s\in[0,T]$ since these processes are c\`adl\`ag.

Below we denote by $C$ a generic constant  depending only on  $T,m,K_1,K_2$, whose precise value may possibly change
 at each occurrence.  Recalling the polynomial growth condition on $f$ and $g$ in {\bf (H2)}, it follows from  Lemma \ref{bsdepenlizaedandcontrol} that
\beqs
|Y_s^n| &\le&  C \mathop{\rm ess\, sup}_{\nu\in\calv^n(\D_\lambda)}  \E^\nu \Big[1+ \sup_{r\in [0,T]}|X_r|^m  \big| \calf_s\Big], \qquad s \in [0,T], \; \P-a.s.
\enqs
Next we note that standard estimates on the stochastic equation  \eqref{defxi} satisfied by $X$, based on the
Lipschitz and linear growth conditions in {\bf (H1)},  lead for every $s \in [0,T]$  to the inequality
\beqs
 \E^\nu \Big[ \sup_{r\in [s,T]}|X_r|^m \big| \calf_s\Big] &\le & C \Big( 1+\sup_{r \in [0,s]} |X_r|^m \Big),
\qquad\P-a.s.
\enqs
which can be viewed as a conditional form of the estimate  \eqref{solusp}, and where the constant $C$ can be chosen to be the same for every $\nu\in \calv(\D_\lambda)$.
It follows that $|Y_s^n| \le C\big(1+ \sup_{r\in [0,s]}|X_r|^m\big)$,  $\P$-a.s.  and the required conclusion follows immediately.
\qed

\begin{Remark}
{\rm  In \cite{khph}, uniform estimates  for $(Y^n)$ in  $\cals^2$ are obtained in the general case where the generator $f$ may also depend on $Y^n,Z^n$, but under the assumption that there exists a
solution to \eqref{BSDEnonpos}. Here, in our specific control case (which leads to the explicit  formula \eqref{expliYn}), one derives directly from Lemma \ref{estimatesonbsdepenlizaed} and \eqref{solusp}  a uniform estimate for $(Y^n)$ in $\cals^2$:
\beqs
\E\Big[  \sup_{s\in [0,T]}|Y_s^n|^2 \Big] & \leq &  C\big( 1 +  \|x \|_\infty^{2m} \big).
\enqs
\qed
}
\end{Remark}

\vspace{2mm}

We are now in a position to complete the proof of Theorem \ref{main2}.

\vspace{3mm}

\noindent {\bf End of the proof of Theorem \ref{main2}}.

\vspace{1mm}

\noindent $\bullet$  Lemma \ref{estimatesonbsdepenlizaed} corresponds to the statements of Lemma 3.1 and 3.2
in \cite{khph}. Following the same arguments used there, it can be proved that  the sequence $(Y^n,Z^n,U^n,n\int_0^\cdot\int_A U^n(a)^+ \,\lambda(da)\,ds)$
converges, in the sense specified in \cite{khph}, Theorem 3.1, to
the required minimal solution  $(Y,Z,$ $U,K)$ to the BSDE  \eqref{BSDEnonpos}.  Thus, for $s\in[0,T]$,  the value
$Y_s$ $=$ $\lim_{n\rightarrow\infty}\uparrow Y_s^n$  inherits from $Y_s^n$ the property of being deterministic, and not dependent on the choice of the set-up $\D_\lambda$.

\vspace{1mm}

\noindent $\bullet$   Since $Y_s^n,X_s,I_s$ are deterministic for $s\in [0,t]$,  it follows from Lemma \ref{bsdepenlizaedandcontrol} that
\beqs
Y_t^n &=& \sup_{\nu\in\calv^n(\D_\lambda)}  \E^\nu \Big[ \int_t^T f_s(X,I_s) \,ds +g(X) \Big].
\enqs
By sending $n$ to infinity, this shows that
\beq \label{vYn}
Y_t \; = \; \lim_{n\rightarrow\infty} \uparrow Y_t^n &=& \sup_{\nu\in\calv(\D_\lambda)}  \E^\nu \Big[ \int_t^T f_s(X,I_s) \,ds +g(X) \Big] \; =  \; v_1^*(t,x,a,\lambda,\D_\lambda).
\enq
Indeed, since  $\calv^n(\D_\lambda)$ $\subset$ $\calv(\D_\lambda)$, we have $\lim_{n\rightarrow\infty} \uparrow Y_t^n$ $\leq$ $v_1^*(t,x,a,\lambda,\D_\lambda)$.
To prove the opposite inequality, take any $\nu\in \calv(\D_\lambda)$ and define $\nu^n_t(a)$ $=$ $\nu_t(a)\wedge n$. By similar arguments to those used in the proof
of equality \eqref{convgainsapprox}, one can show  that for $n\to\infty$
\beqs
 \E^{\nu^n} \Big[ \int_t^Tf_s(X,I_s)\,ds+g(X)\Big]  &\longrightarrow & \E^{\nu} \Big[ \int_s^Tf_s(X,I_s)\,ds+g(X)\Big],
\enqs
which implies that $\lim_{n\to\infty }Y_t^n$ $\ge$ $v_1^*(t,x,a,\lambda, \D_\lambda)$ since $\nu$ is arbitrary in $\calv(\D_\lambda)$, and thus the equality \eqref{vYn}.
This shows in particular that $v_1^*$ does not depend on the choice of $\D_\lambda$, and so is equal to $v_2^*$.  Moreover,  by combining with the result of Theorem \ref{main}, we obtain that  $Y_t$  represents also the  value $v(t,x)$ of the original control problem.
\qed

\appendix

\section{Appendix: some facts on random measures and their compensators}

\setcounter{equation}{0} \setcounter{Assumption}{0}
\setcounter{Theorem}{0} \setcounter{Proposition}{0}
\setcounter{Corollary}{0} \setcounter{Lemma}{0}
\setcounter{Definition}{0} \setcounter{Remark}{0}

We first  recall the following fact, used in the previous sections.

\begin{Proposition}\label{sumofindepmeasures}
Let $(\Omega , \calf ,\P)$ be a probability space endowed with two filtrations $(\calf _t)$,
$(\calk _t)$, such that $\calf _\infty$ and $ \calk _\infty$ are independent.  Let $(T_n,\alpha_n)_{n\ge 1}$ be a
marked point process with respect to $(\calf _t)$ and  $(S_n,\beta_n)_{n\ge 1}$  a
marked point process with respect to $(\calk _t)$.
Denote  $\mu=\sum_{n\ge 1}\delta_{ (T_n,\alpha_n)}$, $\kappa =\sum_{n\ge 1}\delta_{ (S_n,\beta_n)}$  the associated random measures, and $\tilde\mu(dt,da)$,
 $\tilde\kappa(dt,da)$ the respective $(\calf _t)$- and $(\calk _t)$-compensators. Finally assume that the processes have no common jumps.

Then the random measure $\mu+\kappa$ admits
$\tilde\mu(da\,dt)+\tilde\kappa(da\,dt)$ as a
$(\calf _t\vee \calk _t)$-compensator.
 \end{Proposition}

\vspace{2mm}

The rest of this section is devoted to a technical result, stated below as Lemma \ref{MPPperturbed}, that was
used in the proof of Lemma \ref{compensatorofmodifiedppmbis}.

We recall that in our paper we require $A$ to be a Lusin space,
%i.e. a topological space homeomorphic to a Borel subset of a compact metric space (or, equivalently, to a a Borel subset of a separable  complete metric space)
and we denote by $\calb(A)$ its Borel $\sigma$-algebra.

\begin{Lemma}\label{transkernel}
Let $q(b,da)$ be a transition kernel on the Lusin space $A$.
Then there exists a function $q:A\times [0,1]\to A$, measurable with respect to $\calb(A)\otimes \calb([0,1])$ and $\calb(A)$, such that for every
$b\in A$, the measure $B \mapsto q(b,B)$ ($B\in \calb(A)$) is the image of the Lebesgue measure  on $[0,1]$ under the mapping
$u \mapsto q(b,u)$; equivalently,
$$
\int_A k(a)\,q(b,da) \; = \; \int_0^1k( q(b,u))\,du,
$$
for every nonnegative measurable function $k$ on $A$.
\end{Lemma}
{\bf Proof.} When $A$ is a separable complete metric space
(in particular, when $A$ is the unit interval $[0,1]$) the result  is known and follows from a construction of Skorohod:
see e.g. \cite{za}, Theorem 3.1.1. The general case reduces to this one, since it is known that any Lusin  space is either finite or countable (with the discrete topology)
or isomorphic, as a measurable space, to the interval  $[0,1]$: see e.g.
\cite{besh}, Corollary 7.16.1.
\qed

\vspace{2mm}

From the lemma,  it follows that if  $U$ is a random variable defined on some probability
space and having uniform law on $[0,1]$ then, for fixed $b\in A$, the random variable $q(b,U)$ has law $q(b,da)$ on $A$. The use of the same letter $q$ should
not be a source of confusion.

Now let $(\Omega , \calf ,\P)$ be a probability space with a filtration $(\calf _t)$,
let $(T_n,\alpha_n)_{n\ge 1}$ be a marked point process in $A$,
with respect to $(\calf _t)$, and let $\mu=\sum_{n\ge1}\delta_{(T_n,\alpha_n)}$
the corresponding random measure.
Assume that  we can find
sequences $(U_n)_{n\ge1}$, $(S_n)_{n\ge1}$ of
real random variables defined on $\Omega$ and satisfying the following conditions:
\begin{enumerate}
\item every $ U_n $ is uniformly  distributed on $[0,1]$;
\item  every   $ S_n$ admits
 a density (denoted $f_n(t)$)  with respect to the Lebesgue measure,
 and we have
 $
 0<S_1<S_2<S_3<\ldots <\infty$;
 \item the random variables $U_n,S_k$ ($n,k\ge1$) are independent, and independent
 of $\calf_\infty$.
 \end{enumerate}
 Define
 $$
R_n=T_n+S_n,\quad \beta_n=q(\alpha_n,U_n),
\qquad \qquad n\ge 1,
$$
and set
$R_0=0$.
Since we assume $S_n<S_{n+1}$ we see that $(R_n,\beta_n)_{n\ge1}$
is a marked point process in $A$. Let $\kappa=\sum_{n\ge1}\delta_{(R_n,\beta_n)}$
denote the corresponding random measure, and $(\calf^\kappa_t)$ the corresponding
natural filtration. Finally set $\calh_t=\calf_t\vee \calf^\kappa_t$.
We wish to compute the compensator
$\tilde\kappa$ of $\kappa$ with respect to the filtration $\H$ $=$ $(\calh_t)$.
We use the notation $F_n(s)=\int_{-\infty}^sf_n(t)dt$
and the convention $\frac{f_n(s)}{1-F_n(s)}=0$ if $F_n(s)=1$.

\begin{Lemma}\label{MPPperturbed}
With the previous assumptions and notations, the compensator of the random measure
$\kappa$ with respect to $(\calh_t)$ is given by the formula
$$
\tilde\kappa(dt,da) \; = \;  \sum_{n\ge 1}  1_{(T_n\vee R_{n-1},R_n]}(t)\,
 q(\alpha_n,da) \frac{f_n(t-T_n)}{1-F_n(t-T_n)}\,dt.
$$
\end{Lemma}
{\bf Proof.}
Let us first check that $\tilde\kappa(dt,da)$, defined by the  formula above, is an $(\calh_t)$-predictable random measure.
We note that $T_n\vee R_{n-1}$ and $R_n$ are $(\calh_t)$-stopping times
and that $\alpha_n$ and $\frac{f_n(t-T_n)}{1-F_n(t-T_n)}$ are $\calf_{T_n}$-measurable and hence
$\calf_{T_n\vee R_{n-1}}$-measurable. It follows that for every $C\in\calb(A)$ the process
$$
1_{(T_n\vee R_{n-1},R_n]}(t)\, q(\alpha_n,C) \frac{f_n(t-T_n)}{1-F_n(t-T_n)}
$$
is $(\calh_t)$-predictable and finally that  $\tilde\kappa(dt,da)$ is an $(\calh_t)$-predictable random measure.

To prove the Lemma we need now to verify that for every  positive $\calp(\H)\otimes \calb(A)$-measurable
 random field $H_t(\omega,a)$ we have
 $$
 \E \Big[\int_0^\infty\int_AH_t(a)\,\kappa(dt\,da) \Big]  \;= \;
 \E \Big[ \int_0^\infty\int_AH_t(a)\,\tilde\kappa(dt\,da)\Big].
 $$
 Since $\calh_t=\calf_t\vee \calf^\kappa_t$, by a monotone class argument
 it is enough to consider $H$ of the form
 $$H_t(\omega,a)=H_t^1(\omega)H_t^2(\omega)k(a),
 $$
 where $H^1$ is a positive $(\calf_t)$-predictable random process,
$H^2$ is a positive $(\calf^\kappa_t)$-predictable random process
 and $k$ is a positive $\cala$-measurable function.
 Since $(\calf^\kappa_t)$ is the natural filtration of $\kappa$,
 by a known result (see e.g. \cite{ja} Lemma (3.3))
 $H^2$ has the following form:
 \begin{eqnarray*}
    H^2_t &=& b_1(t)1_{(0,R_1]}(t)+
 b_2(\beta_1,R_1,t)1_{( R_1,R_2]}(t)+
 b_3(\beta_1,\beta_2 ,R_1,R_2,t)1_{(R_2,R_3]}(t)+ \\
     && \ldots+
 b_n(\beta_1,\ldots,\beta_{n-1},R_1,\ldots, R_{n-1},t)1_{(R_{n-1},R_n]}(t)+
 \ldots,
 \end{eqnarray*}
where each $b_n$ is a positive measurable deterministic function of
$2n-1$ real variables.
Since
$$
\E \Big[ \int_0^\infty\int_AH_t(a)\,\kappa(dt\,da) \Big] \; = \;
\E \Big[  \sum_{n\ge 1} H_{R_n}(\beta_n) \Big]
 $$
to prove the thesis  it is enough to check that for every $n\ge1$ we have the equality
$$
\E \big[ H_{R_n}(\beta_n) \big] \; = \;
 \E \Big[ \int_0^\infty \int_A H_t(a)\,  q(\alpha_n,da) \frac{f_n(t-T_n)}{1-F_n(t-T_n)}\, 1_{T_n\vee R_{n-1}<t\le R_n }\,dt \Big],
 $$
which can also be written
$$
\begin{array}{l}\dis
 \E\big[H^1_{R_n}  b_n(\beta_1,\ldots,\beta_{n-1},R_1,\ldots, R_{n-1},R_n)    k(\beta_n) \big]
 =
 \\\dis
   \E \Big[ \int_0^\infty \!\int_A H_t^1  b_n(\beta_1,\ldots,\beta_{n-1},R_1,\ldots, R_{n-1},t)k(a)
 q(\alpha_n,da) \frac{f_n(t-T_n)}{1-F_n(t-T_n)}\,  1_{T_n\vee R_{n-1}<t\le R_n }\,dt \Big].
\end{array}
$$
We use the notation
$$
K_n(t) \; = \; H_t^1\,
 b_n(\beta_1,\ldots,\beta_{n-1},R_1,\ldots, R_{n-1},t)
 $$
 to reduce the last equality to
\begin{equation}\label{thesisrewritten}
    \E\, [K_n(R_n)
    k(\beta_n)]    =
 \E\int_0^\infty \int_A K_n (t)\,k(a)\,
 q(\alpha_n,da) \frac{f_n(t-T_n)}{1-F_n(t-T_n)}\,
 1_{T_n\vee R_{n-1}<t\le R_n }\,dt.
\end{equation}
By the definition of $R_n$ and $\beta_n$, we have $\E[K_n(R_n)k(\beta_n)]$ $=$  $\E[K_n(T_n+S_n)   k(q(\alpha_n,U_n))] $.
As noted above,   since $U_n$ has uniform law on $(0,1)$,  the random variable
$q(b,U_n)$ has law $q(b,da)$ on $A$,   for any fixed $b\in A$. Recalling that $S_n$ has density $f_n$ and noting that $U_n $,
 $S_n$  and $(\calf_\infty, \beta_1,\ldots,\beta_{n-1},R_1,\ldots, R_{n-1})$ are all  independent we obtain
\begin{equation}\label{S_nindep}
        \E\, [K_n(R_n)k(\beta_n)]
    \; = \;
    \E \Big[ \int_0^\infty\int_A\, K_n(T_n+s)\,    k(a)\,q(\alpha_n,da)\,f_n(s)\,ds\Big].
\end{equation}
 Using  again the independence of $S_n$  and
$(\calf_\infty, \beta_1,\ldots,\beta_{n-1},R_1,\ldots, R_{n-1})$ we also have
\beqs
& &   \E \Big[ \int_0^\infty\int_A\, K_n(T_n+s)\,    k(a)\,q(\alpha_n,da)\,
    \frac{f_n(s)}{1-F_n(s)}\,1_{S_n\ge s}\, ds \Big] \\
 &=&     \E \Big[ \int_0^\infty\int_A\, K_n(T_n+s)\,    k(a)\,q(\alpha_n,da)\,  \frac{f_n(s)}{1-F_n(s)}\,    \P(S_n\ge s)\, ds\Big]
 \enqs
and since  $\P(S_n\ge s)=\int_s^\infty f_n(r)\,dr=1-F_n(s)$, this coincides with the
right-hand side of \eqref{S_nindep}. By a change of variable we arrive at
the equality
\begin{equation}\label{thesisproved}
\begin{array}{lll}\dis
        \E\, [K_n(R_n)k(\beta_n)]
    &=&\dis
    \E \Big[ \int_{T_n}^\infty\int_A\, K_n(t)\,    k(a)\,q(\alpha_n,da)\,
\frac{f_n(t-T_n)}{1-F_n(t-T_n)}\,    1_{S_n\ge t-T_n}\, dt \Big]
    \\
&=&\dis
    \E \Big[ \int_{0}^\infty\int_A\, K_n(t)\,    k(a)\,q(\alpha_n,da)\,
  \frac{f_n(t-T_n)}{1-F_n(t-T_n)}\,   1_{T_n<t\le R_n}\,dt\Big].
    \end{array}
\end{equation}
We finally claim that
$$
f_n(t-T_n)\, 1_{T_n<t\le R_n}  \; = \;
f_n(t-T_n)1_{T_n\vee R_{n-1}<t\le R_n },
\qquad \P\otimes dt -a.s.
$$
If we can prove the claim, we conclude that \eqref{thesisproved} coincides
with  \eqref{thesisrewritten} and the proof will be finished. To prove the claim, we show that  the following integral is zero:
\begin{equation}\label{funzindic}
    \E \Big[ \int_0^\infty f_n(t-T_n)\, |1_{T_n<t\le R_n}-
1_{T_n\vee R_{n-1}<t\le R_n }| \,dt \Big]
 \; = \;
\E \Big[ \int_0^\infty f_n(t-T_n)\,  1_{T_n<t\le T_n\vee R_{n-1}} \,dt\Big].
\end{equation}
Since $T_n\vee R_{n-1}= T_n\vee (T_{n-1}+S_{n-1})
\le T_n\vee (T_{n}+S_{n-1})=T_{n}+S_{n-1}$,  the right-hand side
of \eqref{funzindic} is smaller or equal to
$$
\E \Big[ \int_0^\infty f_n(t-T_n)\,  1_{T_n<t\le T_n+ S_{n-1}} \,dt \Big]
\; = \;
\E \Big[ \int_0^\infty f_n(s)\,  1_{s\le  S_{n-1}} \,ds \Big].
$$
Since $S_{n}$ and $S_{n-1}$ are independent, and $S_{n-1}<S_n$ a.s.,
we finally have
$$
\E \Big[ \int_0^\infty f_n(s)\,  1_{s\le  S_{n-1}} \,ds\Big]  \; = \;
\E \big[ 1_{S_n\le  S_{n-1}} \big] \; = \;  0.
$$
\qed


\small


\begin{thebibliography}{11}

\bibitem{besh}
Bertsekas, D.P. and S.E. Shreve (1978): Stochastic optimal control. The discrete time case. {\it Mathematics in Science and Engineering}, 139. Academic Press.


\bibitem{bou09} Bouchard, B. (2009):  ``A stochastic target formulation for optimal switching problems in finite horizon", {\it Stochastics}, {\bf 81}, 171-197.


\bibitem{cheetal07} Cheridito, P., Soner, M., Touzi, N and N. Victoir (2007): ``Second-order backward stochastic differential equations and fully nonlinear PDEs", {\it Communication in Pure and Applied Mathematics},
{\bf 60}, 1081-1110.


\bibitem{craishlio92} Crandall M., Ishii H. and P.L. Lions (1992)~: ``User's guide to viscosity solutions of second order partial differential equations", {\it Bull. Amer. Math. Soc.}, {\bf 27}, 1-67.


\bibitem{denetal}
Denis, L.,   Hu, M.  and S. Peng (2011):  ``Function spaces and capacity related to a sublinear expectation: application to $G$-Brownian motion paths",
{\it Potential Anal.}, {\bf 34}, 139-161.


\bibitem{ekretal11}  Ekren, I.,   Keller, K,   Touzi, N. and  and J. Zhang (2011): ``On Viscosity Solutions of Path Dependent PDEs", to appear in {\it Annals of Probability}.


\bibitem{elk79} El Karoui, N. (1981): Les aspects probabilistes du contr\^ole stochastique, {\it Ecole d'\'et\'e de probabilit\'es de St-Flour}, vol. 876,  Lect. Notes in Maths, 73-238, Springer.



\bibitem{FGSbook} Fabbri, G., Gozzi, F. and A. \'Swi\c{e}ch:
 Stochastic Optimal Control in Infinite Dimensions: Dynamic Programming and HJB Equations,  book in preparation, draft.


\bibitem{fleson94} Fleming, W. and M. Soner (2006): Controlled Markov processes and viscosity solutions,  Springer, Series SMAP,  vol. 25.


\bibitem{hewayan}
He, S.W., Wang, J.G. and J.A. Yan (1992):  Semimartingale theory and stochastic calculus.
CRC Press, Boca Raton, FL.


\bibitem{ja}
Jacod, J. (1975):  ``Multivariate point processes: predictable projection,
Radon-Nikodym derivatives, representation of martingales",  {\it Z. Wahrscheinlichkeitstheorie und Verw. Gebiete}, {\bf   31}, 235-253.


\bibitem{khaetal10} Kharroubi I., J. Ma, H. Pham, and J. Zhang (2010): ``Backward SDEs with constrained jumps and Quasi-variational inequalities'', {\it Annals of Probability}, {\bf 38}, 794-840.


\bibitem{khph}
Kharroubi, I. and H.  Pham (2012):  ``Feynman-Kac representation for Hamilton-Jacobi-Bellman IPDE",  Preprint  arXiv:1212.2000.


\bibitem{kr} Krylov, N. V. (1980):  Controlled diffusion processes.
{\it Stochastic Modelling and Applied Probability}, {\bf 14}. Springer-Verlag, Berlin.


\bibitem{nut12} Nutz M. (2012): ``A quasi-sure approach to the control of non-markovian stochastic differential equations", to appear in {\it Electronic Journal of Probability.}



\bibitem{parpen92} Pardoux E. and S. Peng (1992): ``Backward stochastic differential equation and quasilinear parabolic partial differential equations", in {\it Stochastic partial differential equations and their applications}, B. Rozovskii and R. Sowers (eds), Lect. Notes in Cont. Inf. Sci., {\bf 176}, 200-217.



\bibitem{pen06} Peng S. (2006):  ``$G$-expectation, $G$-Brownian motion and related stochastic calculus of It\^o type", Proceedings of 2005, Abel symposium, Springer.

\bibitem{pen11} Peng S. (2011): ``Note on viscosity solutions of path-dependent PDE and $G$-martingales, preprint arXiv: 1106.1144v1


\bibitem{pha09} Pham H. (2009): Continuous-time stochastic control and optimization with financial applications, Springer, Series SMAP, vol, 61.



\bibitem{proshi08} Protter P. and K. Shimbo (2008): ``No arbitrage and general semimartingale", in {\it Festschrift for Thomas Kurtz}.


\bibitem{sontouzha11} Soner M., Touzi N., and J. Zhang (2011):
``The wellposedness of second order backward SDEs",  {\it Probability Theory and Related Fields},
{\bf 153}, 149-190.



\bibitem{rowi}
Rogers, L.C.G. and D.  Williams (2000):  Diffusions, Markov processes, and martingales. Vol. 2. It\^o calculus.
{\it Cambridge University Press}, Cambridge.



\bibitem{lita} Tang, S. and X. Li (1994):  ``Necessary conditions for optimal control of stochastic systems with random jumps",
{\it SIAM J. Control  Optim.}, {\bf  32},  no. 5,   1447-1475.


\bibitem{tanzha12}  Tang S. and F. Zhang (2013): ``Path-dependent optimal stochastic control and viscosity solution of associated Bellman equations", preprint arxiv: 1210.2078

\bibitem{yozh} Yong, J. and Zhou, X.Y. (1999).
Stochastic controls. Hamiltonian systems and HJB equations. Applications of Mathematics
43. Springer, New York.

\bibitem{za}
Zabczyk, J. (1996):  Chance and decision. Stochastic control in discrete time. {\it Quaderni Scuola Normale Superiore}, Pisa.



\end{thebibliography}
\end{document}